\documentclass[12pt]{article}
\usepackage{amsfonts,amscd,a4}
\usepackage{color}
\usepackage{graphicx}
\usepackage{theorem}
\newtheorem{theo}{Theorem}[section]
\newtheorem{prop}[theo]{Proposition}
\newtheorem{lemma}[theo]{Lemma}
\newtheorem{coro}[theo]{Corollary}

{\theorembodyfont{\rm}

}
\newcommand{\bA}{{\bf A}}

\newcommand{\cA}{{\mathcal A}}
\newcommand{\cB}{{\mathcal B}}

\newcommand{\cF}{{\mathcal F}}
\newcommand{\cG}{{\mathcal G}}

\newcommand{\cJ}{{\mathcal J}}

\newcommand{\cO}{{\mathcal O}}

\newcommand{\cS}{{\mathcal S}}

\newcommand{\cY}{{\mathcal Y}}

\newcommand{\eA}{{\sf A}}
\newcommand{\eB}{{\sf B}}

\newcommand{\eD}{{\sf D}}

\newcommand{\eO}{{\sf O}}

\newcommand{\eT}{{\sf T}}

\newcommand{\sC}{{\mathbb C}}

\newcommand{\sF}{{\mathbb F}}

\newcommand{\sN}{{\mathbb N}}

\newcommand{\sZ}{{\mathbb Z}}
\newcommand{\qed}{\rule{1ex}{1ex}}
\newcommand{\alg}{\mbox{\rm alg} \,}

\newcommand{\im}{\mbox{\rm im} \,}

\newcommand{\op}{\mbox{\rm Op} \,}

\begin{document}
\title{Spatial discretization of restricted group algebras}
\author{Steffen Roch}
\date{}
\maketitle
\begin{abstract}
We consider spatial discretizations by the finite section method of the restricted group algebra of a finitely generated discrete group, which is represented as a concrete operator algebra via its left-regular representation. Special emphasis is paid to the quasicommutator ideal of the algebra generated by the finite sections sequences and to the stability of sequences in that algebra. For both problems, the sequence of the discrete boundaries plays an essential role. Finally, for commutative groups and for free non-commutative groups, the algebras of the finite sections sequences are shown to be fractal.
\end{abstract}
\section{Introduction} \label{s1}
Approximately finite algebras and quasi-diagonal algebras are examples of  $C^*$-algebras which are distinguished by intrinsic finiteness properties. These properties can be used in principle to approximate the elements of the algebra by finite-dimensional (or discrete) objects and, thus, to discretize the algebra in a sense. In this paper we consider a completely different kind of discretization, called {\em spatial discretization}, the main idea of which is as follows: We represent a given $C^*$-algebra $\cA$ faithfully as an algebra $\eA$ of linear bounded operators on a separable Hilbert space with basis $\{e_i\}_{i \in \sN}$. Then we let $P_n$ stand for the orthogonal projection from $H$ onto the linear span of $e_1, \, \ldots, \, e_n$, associate with each operator $A \in \eA$ the sequence $(P_n A P_n)$ of its finite sections, and consider the $C^*$-algebra $\cS(\eA)$ which is generated by all sequences $(P_n A P_n)$ with $A \in \eA$. There is a natural homomorphism from $\cS(\eA)$ onto $\eA$ which associates with each sequence in $\cS(\eA)$ its strong limit. Thus, the algebra $\eA$ appears as a quotient of $\cS(\eA)$ by the ideal of all sequences tending {\em strongly} to zero.

The idea of spatial discretization has its origins in numerical analysis, where the numerical solution of an operator equation $Au = f$ is a basic problem. Numerical analysis provides a huge arsenal of methods to discretize this equation for several classes of operators. The perhaps simplest (from the conceptual point of view) and most universal (applicable to each operator) method is the finite sections method which replaces the equation $Au=f$ by the sequence of the finite-dimensional linear systems $P_n A P_n u_n = P_n f$, $n = 1, 2, \ldots$. The basic question is if these systems are uniquely solvable for sufficiently large $n$ and if their solutions $u_n$ tend to a solution of $Au = f$. The central aspect of this question is if the operators ($= n \times n$-matrices) $P_n A P_n$ are invertible for sufficiently large $n$ and if the norms of their inverses are uniformly bounded. In this case, the sequence $(P_n A P_n)$ is called {\em stable}.

A Neumann series argument shows that the sequence $(P_n A P_n)$ with $A \in \eA$ is stable if and only if its coset is invertible in the quotient of the algebra $\cS(\eA)$ by the ideal of all sequences which tend to zero {\em in the norm}. This observation due to Kozak brings numerical analysis into the realm of $C^*$-algebras (and conversely). It was soon realized that, for instance, Gelfand theory and its several non-commutative generalizations provide effective tools to study stability problems for the finite sections method for convolution type equations; see \cite{HRS2} for an overview. In the consequence, the algebras $\cS(\eA)$ were examined for several classes of operator algebras $\eA$. The pioneering example was the Toeplitz algebra, $\eT(C)$, generated all Toeplitz operators on $l^2(\sN)$ with continuous generating function. This algebra can be viewed as a faithful representation of the universal $C^*$-algebra generated by one isometry (Coburn's theorem, \cite{Cob1}). The algebra $\cS(\eT(C))$ of the finite sections method is very well understood; for several aspects of finite sections of Toeplitz operators as well as for the rich history of the field see \cite{BGr5,BSi2}. These results were later extended to algebras generated by Toeplitz operators with piecewise continuous (and even "more discontinuous") symbols and to algebras of singular integral operators, see \cite{HRS1}. The algebra $\cS(\eB \eD \eO)$ of the finite sections of band-dominated operators was subject of \cite{Roc10,RRS6} (note that the algebra $\eB \eD \eO$ of the band-dominated operators is a faithful representation of the reduced crossed product algebra $l^\infty(\sZ) \times_{\alpha r} \sZ$), and the algebra $\cS (\eO_N)$ where $\eO_N$ is a concrete representation of the Cuntz algebra $\cO_N$ was considered in \cite{Roc11}.

The present paper is devoted to the spatial discretization of  restricted group algebras $C^*_r(\Gamma)$ where $\Gamma$ is a finitely generated discrete and exact group. Basic properties of group algebras can be found, e.g., in \cite{Bla2,BroO,Dav1}. Restricted group algebras come with a natural representation, the so-called left-regular representation, which makes $C^*_r(\Gamma)$ isomorphic to the algebra ${\sf Sh}(\Gamma)$ of shift operators on $l^2(\Gamma)$. It is this algebra to which spatial discretization is applied in what follows.

The paper is organized as follows. In Section \ref{c61} we provide some preliminaries on spatial discretization of represented  $C^*$-algebras. Section \ref{c62} is devoted to the spatial discretization of ${\sf Sh}(\Gamma)$. For we choose a family $\cY = (Y_n)$ of finite subsets of $\Gamma$ and consider the sequence of the finite sections $P_{Y_n} A P_{Y_n}$ of $A \in {\sf Sh}(\Gamma)$. We show that the algebra $\cS_\cY ({\sf Sh}(\Gamma))$ generated by these sequences splits into the direct sum of ${\sf Sh}(\Gamma)$ and of an ideal which can be characterized as the quasicommutator ideal of the algebra. A main result is that the sequence $(P_{\partial Y_n})$ of the discrete boundaries always belongs to the algebra $\cS_\cY ({\sf Sh}(\Gamma))$, and that this sequence already generates the quasicommutator ideal. This surprising fact has been already observed in other settings, for example for the algebra $\cS(\eT(C))$ of the finite sections method for the Toeplitz operators (a classical result, closely related to the present paper), but also for the algebra $\cS(\eO_N)$ related with Cuntz algebra (see \cite{Roc11}).

In Section \ref{c63} we derive a necessary and sufficient criterion for the stability of sequences in $\cS_\cY ({\sf Sh}(\Gamma))$. The criterion is formulated of terms of limit operators (see \cite{RRSB,Roc10}). It turns out that it is sufficient to consider limit operators with respect to sequences $\eta$ such that each $\eta_n$ belongs to the boundary of some set $Y_{k_n}$, which gives another hint to the exceptional role of the discrete boundaries. In two special settings (commutative groups and free non-commutative groups) we show moreover that one can restrict to the case when $\eta$ is an (inverse) geodesic path, which implies the fractality of the algebra $\cS_\cY ({\sf Sh}(\Gamma))$ for these groups. We will not present the details, but it should be at least mentioned here that one consequence of fractality is the excellent convergence properties of certain spectral quantities. For example, if a sequence $(A_n)$ belongs to a fractal algebra, then the sets of the singular values (the points in the $\epsilon$-pseudospectrum, the points in the numerical range, respectively) of the $A_n$ converge with respect to the Hausdorff metric. For these and other applications of fractality, see \cite{HRS2,Roc0,Roc10,RoS5}.
\section{Spatial discretization} \label{c61}
\subsection{Hilbert spaces and projections} \label{ss61.1}
For a non-empty finite or countable set $X$, let $l^2(X)$ stand for the Hilbert space of all functions $f : X \to \sC$ with
\[
\|x\|^2 := \sum_{x \in X} |f(x)|^2 < \infty.
\]
For $X = \emptyset$, we define $l^2(X)$ as the space $\{0\}$ consisting of the zero element only. For each subset $Y$ of $X$, we consider $l^2(Y)$ as a closed subspace of $l^2(X)$ in a natural way. The orthogonal projection from $l^2(X)$ to $l^2(Y)$ will be denoted by $P_Y$. Thus, $P_X$ and $P_\emptyset$ are the identity and the zero operator, respectively. For $x \in X$, let $\delta_x$ be the function on $X$ which is 1 at $x$ and 0 at all other points. If $X$ is non-empty, then the family $(\delta_x)_{x \in X}$ forms an orthonormal basis of $l^2(X)$, to which we refer as the {\em standard basis}. \index{basis!standard}

For each sequence $(Y_n)_{n \ge 1}$ of subsets of $X$, define its upper and lower limit as
\[
\limsup Y_n := \cap_{k \ge 1} \cup_{n \ge k} Y_n \quad \mbox{and} \quad \liminf Y_n := \cup_{k \ge 1} \cap_{n \ge k} Y_n.
\]
Thus, $\limsup Y_n$ is the set of all $x \in X$ with $x \in Y_n$ for infinitely many $n$, whereas $\liminf Y_n$ contains all $x \in X$ such that $x \in Y_n$ for all but finitely many $n$. A set sequence $(Y_n)$ is said to converge if $\limsup Y_n = \liminf Y_n$. In this case we denote the upper and lower limit by $\lim Y_n$. The following assertions are easy to check.
\begin{prop}
$(a)$ The sequence $(P_{Y_n})$ of projections converges strongly if and only if the set sequence $(Y_n)$ converges. In this case, $\mbox{\rm s-lim} \, P_{Y_n} = P_{\lim Y_n}$. \\[1mm]
$(b)$ The sequence $(P_{Y_n})$ converges strongly to the identity operator if and only if $\liminf Y_n = X$.
\end{prop}
\begin{coro} \label{c160209.4}
$(a)$ If $Y_n \subseteq Y_{n+1}$ for all $n$, then the sequence $(P_{Y_n})$ converges strongly to $P_{\cup_{n \ge 1}  Y_n}$. \\[1mm]
$(b)$ If $Y_m \cap Y_n = \emptyset$ for all $m \neq n$, then the sequence $(P_{Y_n})$ converges strongly to 0.
\end{coro}
\subsection{Algebras of matrix sequences} \label{ss21}
Let $X$ be as before. Given a sequence $\cY := (Y_n)$ of subsets of $X$, let $\cF_\cY$ denote the set of all bounded sequences $\bA = (A_n)$ of operators $A_n : \im P_{Y_n} \to \im P_{Y_n}$. Equipped with the operations
\[
(A_n) + (B_n) := (A_n + B_n), \quad
(A_n) (B_n) := (A_n B_n), \quad (A_n)^* := (A_n^*)
\]
and the norm $\|\bA\|_{\cF_\cY} := \|A_n\|$, the set $\cF_\cY$ becomes a $C^*$-algebra with identity, and the set $\cG_\cY$ of all sequences $(A_n) \in \cF_\cY$ with $\lim \|A_n\| = 0$ forms a closed ideal of $\cF_\cY$. The relevance of the algebra $\cF_\cY$ and its ideal $\cG_\cY$ in our context stems from the fact (following from a simple Neumann series argument) that a sequence $\bA \in \cF_\cY$ is stable if, and only if, the coset $\bA + \cG_\cY$ is invertible in the quotient algebra $\cF_\cY/\cG_\cY$. Thus, every stability problem is equivalent to an invertibility problem in a suitably chosen $C^*$-algebra.

Let further stand $\cF^C_\cY$ for the set of all sequences $\bA = (A_n)$ of operators $A_n : \im P_{Y_n} \to \im P_{Y_n}$ with the property that the sequences $(A_n P_{Y_n})$ and $(A_n^* P_{Y_n})$ converge strongly. By the uniform boundedness principle, the quantity $\sup \|A_n P_{Y_n}\|$ is finite for every sequence $(A_n)$ in $\cF^C_\cY$. Thus, $\cF^C_\cY$ is a closed and symmetric subalgebra of $\cF_Y$ which contains $\cG$. Note that the mapping
\begin{equation} \label{e91.5}
W : \cF^C_\cY \to L(l^2(X)), \quad \bA \mapsto \mbox{s-lim} \,
A_n P_{Y_n}
\end{equation}
is a $^*$-homomorphism.
\subsection{Spatial discretization of represented algebras} \label{ss22}
Let $\eA$ be a $C^*$-subalgebra of $L(l^2(X))$ (i.e., a represented $C^*$-algebra), and let $\cY := \{Y_n\}$ be a sequence of subsets of $X$. Write $D$ for the mapping of {\em spatial} (= finite sections) {\em discretization}, i.e.,
\begin{equation} \label{e91.6}
D : L(l^2(X)) \to \cF_\cY, \quad A \mapsto (P_{Y_n} A P_{Y_n}),
\end{equation}
and let $\cS_\cY (\eA)$ stand for the smallest closed $C^*$-subalgebra of the algebra $\cF_\cY$ which contains all sequences $D(A)$ with $A \in \eA$. Clearly, $\cS_\cY (\eA)$ is contained in $\cF^C_\cY$, and the mapping $W$ in (\ref{e91.5}) induces a $^*$-homomorphism from $\cS_\cY (\eA)$ onto $\eA$. On this level, one cannot say much about the algebra $\cS_\cY (\eA)$. The little one can say will follow from the following simple facts. A proof is in \cite{Roc11}.
\begin{prop} \label{p91.7}
Let $\eA$ and $\cB$ be $C^*$-algebras, $D : \eA \to \cB$ a symmetric linear contraction, and $W : \cB \to \eA$ a $^*$-homomorphism such that $W(D(A)) = A$ for every $A \in \eA$. Then \\[1mm]
$(a)$ $D$ is an isometry, $D(\eA)$ is a closed linear subspace of $\cB$, and $\alg D(\eA)$, the smallest closed subalgebra of $\cB$ which contains $D(\eA)$, splits into the direct sum
\begin{equation} \label{e91.8}
\alg D(\eA) = D(\eA) \oplus (\ker W \cap \alg D(\eA)).
\end{equation}
Moreover, for every $A \in \eA$,
\begin{equation} \label{e91.9}
\|D(A)\| = \min_{K \in \ker W} \|D(A) + K\|.
\end{equation}
$(b)$ If $\cB = \alg D(\eA)$, then $\ker W$ coincides with the
quasicommutator ideal of $\cB$, i.e., with the smallest closed
ideal of $\cB$ which contains all quasicommutators $D(A_1) D(A_2)
- D(A_1 A_2)$ with $A_1, \, A_2 \in \eA$.
\end{prop}
We shall apply this proposition in the following context: $\eA$ is a $C^*$-subalgebra of $L(l^2(X))$, $\cB$ is the algebra $\cS_\cY(\eA)$, $D$ is the restriction of the discretization (\ref{e91.6}) to $\eA$, and $W$ is the restriction of the homomorphism (\ref{e91.5}) to $\cS_\cY(\eA)$. Then Proposition \ref{p91.7} specializes to the following.
\begin{prop} \label{p91.10}
Let $\eA$ be a $C^*$-subalgebra of $L(l^2(X))$. Then the
finite sections discretization $D : \eA \to \cF_\cY$ is an isometry, and $D(\eA)$ is a closed subspace of the algebra $\cS_\cY(\eA)$. This algebra splits into the direct sum
\[
\cS_\cY(\eA) = D(\eA) \oplus (\ker W \cap \cS_\cY(\eA)),
\]
and for every operator $A \in \eA$ one has
\[
\|D(A)\| = \min_{K \in \ker W} \|D(A) + K\|.
\]
Finally, $\ker W \cap \cS_\cY(\eA)$ is equal to the quasicommutator ideal of $\cS_\cY(\eA)$, i.e., to the smallest closed ideal of $\cS _\cY(\eA)$ which contains all sequences
$(P_{Y_n} A_1 P_{Y_n} A_2 P_{Y_n} - P_{Y_n} A_1 A_2 P_{Y_n})$ with operators $A_1, \, A_2 \in \eA$.
\end{prop}
We denote the ideal $\ker W \cap \cS_\cY(\eA)$ by $\cJ(\eA)$. Since the first item in the decomposition $D(\eA) \oplus \cJ(\eA)$ of $\cS_\cY(\eA)$ is isomorphic (as a linear space) to $\eA$, a main part of the description of the algebra $\cS_\cY(\eA)$ is to identify the ideal $\cJ(\eA)$. Here is a first result which describes $\cJ(\eA)$ in terms of generators of $\eA$. Abbreviate $I - P_A =: Q_A$.
\begin{prop} \label{p030309.1}
Let $\eA$ be a $C^*$-subalgebra of $L(l^2(X))$ and let $E$ be a subset of $\eA$ which generates $\eA$ as a Banach algebra, i.e., the smallest closed subalgebra of $\eA$ which contains $E$ is $\eA$. Then, for each $m \ge 2$ and each choice of operators $A_i \in E$, the sequence
\begin{equation} \label{e030309.2}
(P_{Y_n} A_1 Q_{Y_n} A_2 Q_{Y_n} \ldots \, Q_{Y_n} A_m P_{Y_n})_{n \ge 1}
\end{equation}
belongs to $\cJ(\eA)$, and $\cJ(\eA)$ is the smallest closed ideal of $\cS_\cY(\eA)$ which contains all sequences of the form $(\ref{e030309.2})$.
\end{prop}
{\bf Proof.} First we show per induction that all sequences of the form (\ref{e030309.2}) belong to the quasicommutator ideal $\cJ(\eA)$. This is evident for $m=2$:
\[
(P_{Y_n} A_1 Q_{Y_n} A_2 P_{Y_n}) = (P_{Y_n} A_1 A_2 P_{Y_n}) - (P_{Y_n} A_1 P_{Y_n} A_2 P_{Y_n}).
\]
Suppose the assertion is proved for sequences (\ref{e030309.2}) of length less than $m$. Then
\begin{eqnarray*}
\lefteqn{(P_{Y_n} A_1 Q_{Y_n} \ldots \, Q_{Y_n} A_{m-1} Q_{Y_n} A_m P_{Y_n})} \\
& = & (P_{Y_n} A_1 Q_{Y_n} \ldots \, Q_{Y_n} A_{m-1} A_m P_{Y_n}) \\
&& \quad - (P_{Y_n} A_1 Q_{Y_n} \ldots \, Q_{Y_n} A_{m-1} P_{Y_n}) \, (P_{Y_n} A_m P_{Y_n}).
\end{eqnarray*}
The second sequence on the right-hand side of this equality is in $\cJ(\eA)$ by assumption. Write the first sequence as
\begin{eqnarray*}
\lefteqn{(P_{Y_n} A_1 Q_{Y_n} \ldots \, Q_{Y_n} A_{m-2} Q_{Y_n} A_{m-1} A_m P_{Y_n})} \\
& = & (P_{Y_n} A_1 Q_{Y_n} \ldots \, Q_{Y_n} A_{m-2} A_{m-1} A_m P_{Y_n}) \\
&& \quad - (P_{Y_n} A_1 Q_{Y_n} \ldots \, Q_{Y_n} A_{m-2} P_{Y_n}) \, (P_{Y_n} A_{m-1} A_m P_{Y_n}).
\end{eqnarray*}
Again, the second sequence on the right-hand side is in $\cJ(\eA)$. We continue in this way to arrive finally at
\begin{eqnarray*}
(P_{Y_n} A_1 Q_{Y_n} A_2 A_3 \ldots \, A_m P_{Y_n})
& = & (P_{Y_n} A_1 A_2 \ldots \, A_m P_{Y_n}) \\
&& \quad - (P_{Y_n} A_1 P_{Y_n}) \, (P_{Y_n} A_2 A_3 \ldots \, A_m P_{Y_n})
\end{eqnarray*}
which is in $\cJ(\eA)$ by the definition of the quasicommutator ideal.

Conversely, we are going to show that the sequences (\ref{e030309.2}) generate $\cJ(\eA)$ as a closed ideal of
$\cS_\cY(\eA)$. Let $\cJ$ refer to the smallest closed ideal of
$\cS_\cY(\eA)$ which contains all sequences (\ref{e030309.2}). From the first part of this proof we infer that $\cJ \subseteq
\cJ(\eA)$. For the reverse inclusion it is sufficient to show that \[
(P_{Y_n} A Q_{Y_n} B P_{Y_n}) \in \cJ \quad \mbox{for all} \; A, \, B \in \eA.
\]
Since $\cJ$ is a closed linear space, it is sufficient to verify this claim in case $A$ and $B$ are finite products of operators in $E$. Thus, we have to prove that
\begin{equation} \label{e030309.3}
(P_{Y_n} A_1 \ldots A_m Q_{Y_n} B_1 \ldots B_l P_{Y_n}) \in \cJ
\end{equation}
for arbitrary operators $A_i, \, B_j \in E$ and integers $l, \, m \ge 1$. Again we use induction. The assertion is evident in case $m=l=1$. For the general step we write
\begin{eqnarray*}
\lefteqn{(P_{Y_n} A_1 \ldots A_m Q_{Y_n} B_1 \ldots B_l P_{Y_n})} \\
&& = (P_{Y_n} A_1 P_{Y_n}) \, (P_{Y_n} A_2 \ldots A_m Q_{Y_n} B_1 \ldots B_l P_{Y_n}) \\
&& \quad + \; (P_{Y_n} A_1 Q_{Y_n} A_2 \ldots A_m Q_{Y_n} B_1 \ldots B_l P_{Y_n}).
\end{eqnarray*}
The first summand on the right-hand side is a product of a sequence in $\cS_\cY (\eA)$ and a sequence of the form (\ref{e030309.3}), but with less factors. By assumption, this summand is in $\cJ$. The second summand can be again written as a sum by inserting $I = P_{Y_n} + Q_{Y_n}$ after $A_2$. We continue in this way and arrive finally at the sequence $(P_{Y_n} A_1 Q_{Y_n} A_2 Q_{Y_n} \ldots Q_{Y_n} B_{l-1} Q_{Y_n} B_l P_{Y_n})$ which is in $\cJ$ by definition. \hfill \qed
\section[restricted group algebras]{Spatial discretization of restricted group $C^*$-algebras} \label{c62}
\subsection{Left regular representations} \label{ss62.1}
Let $\Gamma$ be a (not necessarily commutative) discrete group. We write the group operation as multiplication and let $e$ stand for the identity element. With $\Gamma$ we associate the Hilbert space $l^2(\Gamma)$ with its canonical basis $(\delta_s)_{s \in \Gamma}$. The {\em left regular representation} \index{representation!left regular} $L : \Gamma \to L(l^2(\Gamma))$ associates with every group element $r$ a unitary operator $L_r$ such that $L_r \delta_s = \delta_{rs}$ for $s \in \Gamma$.

Since $\delta_{rs}(t) = \delta_s(r^{-1}t)$, one has $(L_r u)(t) = u(r^{-1}t)$ for every $u \in l^2(\Gamma)$. Hence, $r \mapsto L_r$ is a group isomorphism. We define ${\sf Sh}(\Gamma)$ as the smallest closed subalgebra of $L(l^2(\Gamma))$ which contains all operators $L_t$ with $t \in \Gamma$. The algebra ${\sf Sh}(\Gamma)$ is $^*$-isomorphic to the {\em restricted group $C^*$-algebra} $C^*_r(\Gamma)$ in a natural way (see Section 2.5 in \cite{BroO}). It can thus be considered as a concrete representation of $C^*_r(\Gamma)$. Note also that the restricted group $C^*$-algebra coincides with the universal group $C^*$-algebra $C^*(\Gamma)$ if the group $\Gamma$ is amenable. For this and further characterizations of amenable groups, see Theorem 2.6.8 in \cite{BroO}.

We have seen above that every {\em restricted} group $C^*$-algebra $C^*_r(\Gamma)$ comes with a canonical faithful representation as the concrete operator algebra ${\sf Sh}(\Gamma)$ on $l^2(\Gamma)$. We will take this representation as the basis for the spatial discretization of $C^*_r(\Gamma)$ by a finite sections method in the following sections.

The existence of a canonical representation is only one reason why we consider spatial discretizations only for restricted group $C^*$-algebras in what follows. Another reason is that universal group $C^*$-algebras sometimes own intrinsic finiteness properties which can be used to approximate their elements by finite dimensional objects, but which are not shared by the associated restricted group $C^*$-algebras. For example, if $\Gamma$ is the free non-commutative group $\sF_2$ of two generators, then the universal group $C^*$-algebra $C^*(\sF_2)$ is known to be quasidiagonal, whereas $C^*_r(\sF_2)$ fails to have this property (see Sections VII.6 and VII.7 in \cite{Dav1}).
\subsection{Discretization of ${\sf Sh}(\Gamma)$}
To discretize the algebra ${\sf Sh}(\Gamma)$ by the finite sections method we choose a sequence $\cY = (Y_n)$ of finite subsets of $\Gamma$ and consider the sequences $(P_{Y_n} A P_{Y_n})$ of the finite sections of $A \in {\sf Sh}(\Gamma)$. Usually we will assume that the set limit $\lim Y_n$ exists and is equal to $\Gamma$, in which case the $P_{Y_n}$ converge strongly to the identity operator, but some of the following results will hold without this assumption. In accordance with earlier notation, let $\cS_\cY ({\sf Sh}(\Gamma))$ stand for the smallest closed $C^*$-subalgebra of the algebra $\cF_\cY$ which contains all sequences $(P_{Y_n} A P_{Y_n})$ with $A \in {\sf Sh}(\Gamma)$. The associated quasicommutator ideal is denoted by $\cJ ({\sf Sh}(\Gamma))$.

In the next section, we shall present some characterizations of $\cJ ({\sf Sh}(\Gamma))$. For we have to introduce some notions of topological type. Note that the standard topology on $\Gamma$ is the discrete one; so every subset of $\Gamma$ is open with respect to this topology.

Let $\Omega$ be a finite subset of $\Gamma$ which contains the identity element $e$ and which generates $\Gamma$ as a semi-group, i.e., if $\Omega_n$ denotes the set of all words of length at most $n$ with letters in $\Omega$, then $\cup_{n \ge 0} \Omega_n = \Gamma$. By convention, $\Omega_0 := \{e\}$. Note also that the sequence $(\Omega_n)$ is increasing; so the operators $P_{\Omega_n}$ can play the role of the finite sections projections $P_{Y_n}$, and in fact we will obtain some of the subsequent results exactly for this sequence.

With respect to $\Omega$, we define the following "algebro-topological" notions. Let $A \subseteq \Gamma$. A point $a  \in A$ is called an {\em $\Omega$-inner} point of $A$ if $\Omega a := \{ \omega a : \omega \in \Omega \} \subseteq A$. The set $\mbox{int}_\Omega A$ of all $\Omega$-inner points of $A$ is called the {\em $\Omega$-interior} of $A$, and the set $\partial_\Omega A := A \setminus \mbox{int}_\Omega A$ is the {\em $\Omega$-boundary} of $A$. Note that by this definition, the $\Omega$-boundary of a set is always a part of that set. (In this point, the present definition of a boundary differs from other definitions used in the literature, see, e.g., \cite{Ada1}.)

One easily checks that the $\Omega$-interior and the $\Omega$-boundary of a set are invariant with respect to multiplication from the right-hand side:
\[
(\mbox{int}_\Omega A) s = \mbox{int}_\Omega (As) \quad \mbox{and} \quad (\partial_\Omega A) s = \partial_\Omega (As)
\]
for $s \in \Gamma$. One also has
\begin{equation} \label{e030309.5}
\Omega_{n-1} \subseteq \mbox{int}_\Omega \Omega_n \subseteq \Omega_n \quad \mbox{for each} \; n \ge 1,
\end{equation}
whence
\begin{equation} \label{e030309.6}
\partial_\Omega \Omega_n \subseteq \Omega_n \setminus \Omega_{n-1} \quad \mbox{for each} \; n \ge 1.
\end{equation}
In many concrete settings, one has equality in (\ref{e030309.6}).
\subsection{The structure of the quasicommutator ideal}
Let $\Omega$ and $\cY := (Y_n)$ be as in the previous section. We will derive two results on the structure of the quasicommutator ideal $\cJ ({\sf Sh}(\Gamma))$.
\begin{theo} \label{t030309.8}
$\cJ ({\sf Sh}(\Gamma))$ is the smallest closed ideal of $\cS_\cY ({\sf Sh}(\Gamma))$ which contains all sequences
\begin{equation} \label{e030309.9}
(P_{Y_n} L_{\omega^{-1}} Q_{Y_n} L_\omega P_{Y_n})_{n \ge 1} \quad \mbox{with} \quad \omega \in \Omega.
\end{equation}
\end{theo}
{\bf Proof.} First note that, for arbitrary $\omega \in \Omega$ and $A \subseteq \Gamma$,
\begin{equation} \label{e030309.7}
Q_A L_\omega P_A = Q_A L_\omega P_A L_{\omega^{-1}} Q_A L_\omega P_A.
\end{equation}
Indeed,
\begin{eqnarray*}
Q_A L_\omega P_A & = & Q_A L_\omega P_A L_{\omega^{-1}} L_\omega P_A \\
& = & Q_A L_\omega P_A L_{\omega^{-1}} Q_A L_\omega P_A
+ Q_A L_\omega P_A L_{\omega^{-1}} P_A L_\omega P_A.
\end{eqnarray*}
The second summand on the right-hand side vanishes since $L_\omega P_A L_{\omega^{-1}} = P_{\omega A}$ commutes with $P_A$.

Let now, for a moment, $\cJ$ denote the smallest closed ideal of $\cS_\cY ({\sf Sh}(\Gamma))$ which contains all sequences (\ref{e030309.9}). Clearly, $\cJ \subseteq \cJ ({\sf Sh}(\Gamma))$. The reverse implication will follow via
Proposition \ref{p030309.1} once we have shown that each sequence
\begin{equation} \label{e030309.10}
(P_{Y_n} L_{\omega_1} Q_{Y_n} L_{\omega_2} Q_{Y_n} \ldots
Q_{Y_n} L_{\omega_m} P_{Y_n})_{n \ge 1}
\end{equation}
with $m \ge 2$ and $\omega_i \in \Omega$ belongs to $\cJ$.
Write the sequence (\ref{e030309.10}) as $(A_n Q_{Y_n} L_{\omega_m} P_{Y_n})$. By (\ref{e030309.7}),
\[
(A_n Q_{Y_n} L_{\omega_m} P_{Y_n}) = (A_n Q_{Y_n} L_{\omega_m} P_{Y_n}) \, (P_{Y_n} L_{\omega^{-1}} Q_{Y_n} L_{\omega_m} P_{Y_n}).
\]
Since the sequence (\ref{e030309.10}) belongs to $\cS_\cY ({\sf Sh}(\Gamma))$ and $\cJ$ is an ideal of that algebra, the sequence (\ref{e030309.10}) is in $\cJ$. \hfill \qed
\begin{lemma} \label{l030309.11}
Let $A \subseteq \Gamma$. Then $\cap_{\omega \in \Omega} (A \cap \omega^{-1} A) = \mbox{\rm int}_\Omega A$.
\end{lemma}
{\bf Proof.} Let $a \in \mbox{int}_\Omega A$. Then, for each  $\omega \in \Omega$,
\[
a = \omega^{-1} \omega a \in \omega^{-1} \Omega a \subseteq \omega^{-1} A \subseteq A \cap \omega^{-1} A,
\]
whence the inclusion $\supseteq$. For the reverse inclusion, let $a \in A \setminus \mbox{int}_\Omega A = \partial_\Omega A$. By definition of the $\Omega$-boundary, there is an $\omega_0 \in \Omega$ such that $\omega_0 a \not\in A$. Hence, $a \not\in A \cap \omega_0^{-1} A$, which implies $a \not\in \cap_{\omega \in \Omega} (A \cap \omega^{-1} A)$. \hfill \qed
\begin{lemma} \label{l030309.13}
Let $\eA$ be a subalgebra of $L(l^2(\Gamma))$ and $A \subseteq \Gamma$. If the operators $P_A L_{\omega^{-1}} P_A L_\omega P_A$ belong to $\eA$ for each $\omega \in \Omega$, then the operators $P_A$, $P_{\mbox{\scriptsize \rm int}_\Omega A}$ and $P_{\partial_\Omega A}$ belong to $\eA$, too.
\end{lemma}
{\bf Proof.} Since $e \in \Omega$, the assertion is evident for $P_A$. Further we have
\[
P_A L_{\omega^{-1}} P_A L_\omega P_A = P_A P_{\omega^{-1} A} = P_{A \cap \omega^{-1} A} \in \eA
\]
for each $\omega \in \Omega$. Since $\eA$ is an algebra, this implies
\[
\prod_{\omega \in \Omega} P_{A \cap \omega^{-1} A} = P_{\cap_{\omega \in \Omega} (A \cap \omega^{-1} A)} \in \eA.
\]
By Lemma \ref{l030309.11}, this is the assertion for $P_{\mbox{\scriptsize int}_\Omega A}$. The assertion for $P_{\partial_\Omega A}$ follows since $P_A = P_{\mbox{\scriptsize int}_\Omega A} + P_{\partial_\Omega A}$. \hfill \qed \\[3mm]
We call $(P_{\partial_\Omega Y_n})_{n \ge 1}$ the {\em sequence of the discrete boundaries} of the finite section method with respect to $(Y_n)$. Note that the assumptions in the following theorem are satisfied if $Y_n = \Omega_n$ due to (\ref{e030309.5}).
\begin{theo} \label{t030309.14}
Assume that $Y_{n-1} \subseteq \mbox{\rm int}_\Omega Y_n \subseteq Y_n$ for all $n \ge 2$ and that $\lim Y_n = \Gamma$. Then the sequence $(P_{\partial_\Omega Y_n})_{n \ge 1}$ of the discrete boundaries belongs to the algebra $\cS_\cY ({\sf Sh}(\Gamma))$, and the quasicommutator ideal is generated by this sequence, i.e., $\cJ ({\sf Sh}(\Gamma))$ is the smallest closed ideal of $\cS_\cY ({\sf Sh}(\Gamma))$ which contains $(P_{\partial_\Omega Y_n})_{n \ge 1}$.
\end{theo}
{\bf Proof.} By definition, the sequence $(P_{Y_n} L_{\omega^{-1}} P_{Y_n} L_\omega P_{Y_n})_{n \ge 1}$ is in $\cS_\cY ({\sf Sh}(\Gamma))$ for each $\omega \in \Omega$. From Lemma \ref{l030309.13} we then conclude that the sequence $(P_{\partial_\Omega Y_n})$ is in $\cS_\cY ({\sf Sh}(\Gamma))$, too. That this sequence is even in the quasicommutator ideal, is a consequence of the assumptions. Indeed, from $Y_{n-1} \subseteq \mbox{\rm int}_\Omega Y_n \subseteq Y_n$ we conclude that
\[
\lim_{n \to \infty} \mbox{\rm int}_\Omega Y_n = \lim_{n \to \infty} Y_n = \Gamma
\]
whence $\mbox{s-lim} P_{\partial_\Omega Y_n} = 0$. By Proposition \ref{p91.10}, this implies $(P_{\partial_\Omega Y_n}) \in \cJ ({\sf Sh}(\Gamma))$.

It remains to show that the sequence $(P_{\partial_\Omega Y_n})$ generates $\cJ ({\sf Sh}(\Gamma))$. Let $\cJ$ denote the smallest closed ideal of $\cS_\cY ({\sf Sh}(\Gamma))$ which contains the sequence $(P_{\partial_\Omega Y_n})$. By what we have just seen, $\cJ \subseteq \cJ ({\sf Sh}(\Gamma))$. The reverse inclusion will follow from Theorem \ref{t030309.8} once we have shown that
\begin{equation} \label{e030309.15}
(P_{Y_n} L_{\omega^{-1}} Q_{Y_n} L_\omega P_{Y_n})_{n \ge 1} \in \cJ \quad \mbox{for each} \quad \omega \in \Omega.
\end{equation}
Note that
\[
P_{Y_n} L_{\omega^{-1}} Q_{Y_n} L_\omega P_{Y_n} = P_{Y_n} - P_{Y_n} L_{\omega^{-1}} P_{Y_n} L_\omega P_{Y_n} = P_{Y_n \setminus (Y_n \cap \omega^{-1} Y_n)}.
\]
From Lemma \ref{l030309.11} we know that $\mbox{int}_\Omega Y_n \subseteq Y_n \cap \omega^{-1} Y_n$. Hence,
\[
Y_n \setminus (Y_n \cap \omega^{-1} Y_n) \subseteq Y_n \setminus \mbox{int}_\Omega Y_n = \partial_\Omega Y_n
\]
which implies that
\[
P_{Y_n} L_{\omega^{-1}} Q_{Y_n} L_\omega P_{Y_n} =
P_{Y_n \setminus (Y_n \cap \omega^{-1} Y_n)} = P_{Y_n \setminus (Y_n \cap \omega^{-1} Y_n)} \, P_{\partial_\Omega Y_n}.
\]
This verifies (\ref{e030309.15}) and finishes the proof of the theorem. \hfill \qed
\section{Stability} \label{c63}
In this section, we are going to study the stability of sequences in $\cS_\cY ({\sf Sh}(\Gamma))$ via the limit operators method. The key ingredients are the facts that the stability of a sequence $\bA$ in that algebra is equivalent to the Fredholmness of a certain associated operator and that the Fredholmness of that operator can be studied by means of its limit operators due to a result of Roe.
\subsection{Fredholmness vs. stability}
Let $\cY := (Y_n)$ be a sequence of finite subsets of $\Gamma$. A sequence $(v_n) \subseteq \Gamma$ is called an {\em inflating} sequence for $\cY$ if
$Y_m v_m^{-1} \cap Y_n v_n^{-1} = \emptyset$ for $m \neq n$. The existence of inflating sequences is a consequence of the following lemma.
\begin{lemma}
Let $A, \, B \subset \Gamma$ be finite and $V \subset \Gamma$ be infinite. Then there is a $v \in V$ such that $A \cap Bv^{-1} = \emptyset$.
\end{lemma}
Indeed, let $A \cap Bv^{-1} \neq \emptyset$ for every $v \in V$. Then, for each $v \in V$, there is a $b_v \in B$ such that $b_v v^{-1} =: a_v \in A$. Thus, $v = b_v a_v^{-1}$. But since $A$ and $B$ are finite, there are only finitely many products $b_v a_v^{-1}$. Hence $V$ is finite, a contradiction. \hfill \qed
\begin{coro}
Let $\cY = (Y_n)$ be a sequence of finite subsets of $\Gamma$ and $V$ an infinite subset of $\Gamma$. Then there is an inflating sequence for $\cY$ in $V$.
\end{coro}
{\bf Proof.} Let $v_1 \in V$. Then $Y_1 v_1^{-1}$ is finite. By the lemma, there is a $v_2 \in V$ such that $Y_1 v_1^{-1} \cap Y_2 v_2^{-1} = \emptyset$. Further, since $Y_1 v_1^{-1} \cup Y_2 v_2^{-1}$ is finite, there is a $v_3 \in V$ such that
\[
\left( Y_1 v_1^{-1} \cup Y_2 v_2^{-1} \right) \cap Y_3 v_3^{-1} = \emptyset.
\]
We proceed in this way to find the desired inflating sequence. \hfill \qed \\[3mm]
In what follows let $\cY$ as above and choose and fix an inflating sequence $(v_n)$ for $\cY$. Further set
\begin{equation} \label{e180209.2}
\Gamma^\prime := \Gamma \setminus \cup_{n=1}^\infty Y_n v_n^{-1}.
\end{equation}
For $s \in \Gamma$, let $R_s : l^2(\Gamma) \to l^2(\Gamma)$ refer to the operator $(R_s f)(t) := f(ts)$. Evidently, the mapping $R : s \mapsto R_s$ is a group isomorphism from $\Gamma$ into the group of the unitary operators on $l^2(\Gamma)$. Moreover, $R_s L_t = L_t R_s$ for $s, \, t \in \Gamma$. The proof of the following theorem is adapted from \cite{Roc10}.
\begin{theo} \label{t160209.6}
Let $\bA = (A_n) \in \cF_\cY$. Then \\[1mm]
$(a)$ the series
\begin{equation} \label{e160209.0}
\sum_{n=1}^\infty R_{v_n} A_n R_{v_n}^{-1}
\end{equation}
converges strongly on $l^2(\Gamma)$. The sum of this series is denoted by $\op (\bA)$. \\[1mm]
$(b)$ the sequence $(A_n)$ is stable if and only if the operator $\op (\bA) + P_{\Gamma^\prime}$ is Fredholm on $l^2(\Gamma)$. \\[1mm]
$(c)$ The mapping $\op$ is a continuous homomorphism from $\cF_\cY$ to $L(l^2(\Gamma))$.
\end{theo}
{\bf Proof.} $(a)$ It is convenient to identify the operator $A_n$ acting on $\im P_{Y_n}$ with the operator $P_{Y_n} A_n P_{Y_n}$ acting on all of $l^2(\Gamma)$. Since $R_{v_n} P_{Y_n} R_{v_n}^{-1} = P_{Y_n v_n^{-1}}$, one can then identify the operator
\[
R_{v_n} A_n R_{v_n}^{-1} : \im P_{Y_n v_n^{-1}} \to \im P_{Y_n v_n^{-1}}
\]
with the operator $P_{Y_n v_n^{-1}} R_{v_n} A_n R_{v_n}^{-1} P_{Y_n v_n^{-1}}$ on $l^2(\Gamma)$. Thus, for $x \in l^2(\Gamma)$, the inflating property ensures that the vectors $R_{v_n} A_n R_{v_n}^{-1} x$ form an orthogonal system in $l^2(\Gamma)$. Consequently, the series $\sum_{n=1}^\infty R_{v_n} A_n R_{v_n}^{-1} x$ converges if and only if the series
\begin{equation} \label{160209.1}
\sum_{n=1}^\infty \|R_{v_n} A_n R_{v_n}^{-1} x\|^2
\end{equation}
converges. Set $M := \sup \|A_n\|$. Employing the orthogonality of the vectors $P_{Y_n v_n^{-1}} x$, we get
\[
\sum_{n=1}^\infty \|R_{v_n} A_n R_{v_n}^{-1} x\|^2 \le M^2 \sum_{n=1}^\infty \|P_{Y_n v_n^{-1}} x\|^2 \le M^2 \|x\|^2.
\]
Thus, the series (\ref{160209.1}) converges for every $x$, whence assertion $(a)$. \\[3mm]
$(b)$ Let $\bA = (A_n)$ be a stable sequence, i.e., there is an $n_0 \in \sN$ such that the operators $A_n : \im P_{Y_n} \to \im P_{Y_n}$ are invertible for $n \ge n_0$ and that the norms of their inverses are uniformly bounded. Then the operator
\[
B := \sum_{n=1}^{n_0 - 1} P_{Y_n v_n^{-1}} + \sum_{n = n_0}^\infty R_{v_n} A_n R_{v_n}^{-1} + P_{\Gamma^\prime}
\]
is invertible with inverse
\[
B^{-1} = \sum_{n=1}^{n_0 - 1} P_{Y_n v_n^{-1}} + \sum_{n = n_0}^\infty R_{v_n} A_n^{-1} R_{v_n}^{-1} + P_{\Gamma^\prime}.
\]
Since $\op (\bA) + P_{\Gamma^\prime}$ is a compact perturbation of $B$, $\op (\bA) + P_{\Gamma^\prime}$ is a Fredholm operator (with Fredholm index 0).

Let, conversely, $\op (\bA) + P_{\Gamma^\prime}$ be a Fredholm operator. Then there are an operator $B \in L(l^2(\Gamma))$ and a compact operator $K$ on $l^2(\Gamma)$ such that
\[
B \cdot (\op (\bA) + P_{\Gamma^\prime}) = I + K.
\]
Since the projections $P_{Y_n v_n^{-1}}$ commute with $\op (\bA)$, and since $Y_n v_n^{-1} \cap \Gamma^\prime = \emptyset$, we find
\begin{eqnarray*}
P_{Y_n v_n^{-1}} B P_{Y_n v_n^{-1}} \cdot R_{v_n} A_n R_{v_n}^{-1} & = & P_{Y_n v_n^{-1}} B P_{Y_n v_n^{-1}} \cdot P_{Y_n v_n^{-1}} \op (\bA) P_{Y_n v_n^{-1}} \\
& = & P_{Y_n v_n^{-1}} B \op (\bA) P_{Y_n v_n^{-1}} \\
& = & P_{Y_n v_n^{-1}} B (\op (\bA) + P_{\Gamma^\prime}) P_{Y_n v_n^{-1}} \\
& = & P_{Y_n v_n^{-1}} + P_{Y_n v_n^{-1}} K P_{Y_n v_n^{-1}},
\end{eqnarray*}
whence
\begin{equation} \label{e160209.2}
R_{v_n}^{-1} P_{Y_n v_n^{-1}} B P_{Y_n v_n^{-1}} R_{v_n} \cdot A_n = P_{Y_n} + R_{v_n}^{-1} P_{Y_n v_n^{-1}} K P_{Y_n v_n^{-1}} R_{v_n}.
\end{equation}
Since $P_{Y_n v_n^{-1}} \to 0$ strongly by the inflating property and by Corollary \ref{c160209.4} $(b)$ and since $K$ is compact and $\|R_{v_n}\| = 1$, we further conclude
\[
\|R_{v_n}^{-1} P_{Y_n v_n^{-1}} K P_{Y_n v_n^{-1}} R_{v_n}\| \to 0 \quad \mbox{as} \; n \to \infty.
\]
Hence, the operators on the right-hand side of (\ref{e160209.2}) (considered as acting on $\im P_{Y_n}$) are invertible for $n$ large enough, and the norms of their inverses are uniformly bounded with respect to $n$. This implies the uniform boundedness of the operators
\[
B_n := \left( P_{Y_n} + R_{v_n}^{-1} P_{Y_n v_n^{-1}} K P_{Y_n v_n^{-1}} R_{v_n} \right)^{-1} R_{v_n}^{-1} P_{Y_n v_n^{-1}} B P_{Y_n v_n^{-1}} R_{v_n},
\]
also considered as acting on $\im P_{Y_n}$. Since $B_n A_n = P_{Y_n}$ for all sufficiently large $n$ and the $A_n$ act on a finite-dimensional space, the stability of the sequence $(A_n)$ follows. Assertion $(c)$ is an immediate consequence of the inflating property. \hfill \qed
\subsection{Band-dominated operators}
Theorem \ref{t160209.6} translates the stability problem for a bounded sequence of finite-rank operators into a Fredholm problem for an associated operator. In case of the finite sections sequence of an operator in ${\sf Sh}(\Gamma)$, the associated operator is a band-dominated operator in the sense defined below. Since there is an effective criterion to verify the Fredholm property (which we will recall in the subsequent section) of band-dominated operators, this observation offers a way to study the stability of the finite sections method for operators in ${\sf Sh}(\Gamma)$.

Consider functions $k \in l^\infty (\Gamma \times \Gamma)$ with the property that there is a finite subset $\Gamma_0$ of $\Gamma$ such that $k(t, \, s) = 0$ whenever $t s^{-1} \not\in \Gamma_0$. Then
\begin{equation} \label{1.1groups}
(Au)(t) := \sum_{s \in \Gamma} k (t, \, s) \, u(s), \qquad t \in \Gamma,
\end{equation}
defines a linear operator $A$ on the linear space of all
functions $u : \Gamma \to \sC$, since the occurring series is finite for every $t \in G$. We call operators of this form
{\em band operators} and the set $\Gamma_0$ a {\em band-width} of $A$. It is not hard to see that the band operators form a symmetric algebra of bounded operators on $l^2(\Gamma)$. Operators in the norm closure of that algebra are called {\em band-dominated} operators. Thus, the band-dominated operators form a $C^*$-subalgebra ${\sf BDO}(\Gamma)$ of $L(l^2(\Gamma))$.

It turns out that band operators on $\Gamma$ are constituted by two kinds of "elementary" band operators: the unitary operators $L_t$ of left shift by $t \in \Gamma$, and the operators $bI$ of multiplication by a function $b \in l^\infty (G)$,
\[
bI : l^2(\Gamma) \to l^2(\Gamma), \quad (bu)(s) = b(s) \, u(s).
\]
\begin{prop}
A operator in $L(l^2(\Gamma))$ is a band operator if and only if it can be written as a finite sum $\sum b_i L_{t_i}$ where $b_i \in l^\infty(\Gamma)$ and $t_i \in \Gamma$.
\end{prop}
{\bf Proof.} Let $A$ be an operator of the form (\ref{1.1groups}) and let $\Gamma_0 := \{t_1, \, t_2, \, \ldots, \, t_r \}$ be a finite subset of $\Gamma$ such that $k(t, \, s) = 0$ if $t s^{-1} \not\in \Gamma_0$ or, equivalently, if $s$ is not of the form
$t_i^{-1} t$ for some $i$. Thus,
\[
(Au)(t) = \sum_{i = 1}^r k (t, \, t_i^{-1} t) \, u(t_i^{-1} t) \quad \mbox{for all} \; t \in \Gamma.
\]
Set $b_i(t) := k (t, \, t_i^{-1} t)$. The functions $b_i$ are in $l^\infty (\Gamma)$, and one has
\begin{equation} \label{bdo}
A = \sum_{i = 1}^r b_i L_{t_i}.
\end{equation}
Conversely, one easily checks that each operator $L_t$ with $t \in \Gamma$ is a band operator with band width $\{t\}$ and that each operator $bI$ with $b \in l^\infty(\Gamma)$ is a band operator with band width $\{e\}$. Since the band operators form an algebra, each finite sum $\sum b_i L_{t_i}$ is a band operator. \hfill \qed \\[3mm]
It is easy to see that the representation of a band operator on $\Gamma$ in the form (\ref{bdo}) is unique. The functions $b_i$ are called the {\em diagonals} of the operator $A$. In particular, operators in ${\sf Sh}(\Gamma)$ can be considered as band-dominated operators with constant coefficients.

As before, let $\cY := (Y_n)$ be a sequence of finite subsets of $\Gamma$ and $(v_n)$ an associated inflating sequence. Note that the following proposition remains valid if the algebra ${\sf Sh}(\Gamma)$ is replaced by the $C^*$-algebra ${\sf BDO}(\Gamma)$ of all band-dominated operators.
\begin{prop} \label{p170209.1}
Let $\bA = (A_n)$ be a sequence in the finite sections algebra $\cS_\cY ({\sf Sh}(\Gamma))$. Then $\op (\bA)$ is a band-dominated operator.
\end{prop}
{\bf Proof.} First let $A \in {\sf Sh}(\Gamma)$ be a band operator (i.e., $A$ is a linear combination of a finite number of the $L_t$) and let $\Gamma_0$ be a band width of $A$. It is easy to check that then $R_{v_n} P_{Y_n} A P_{Y_n} R_{v_n}^{-1}$ is a band operator with the same band width for every $n$. The inflating property ensures that $\op \left( (P_{Y_n} A P_{Y_n}) \right)$ is a band operator with band width $\Gamma_0$, too. Now Theorem \ref{t160209.6} $(c)$ yields the assertion. \hfill \qed \\[3mm]
To define limit operators, let $h : \sN \to \Gamma$ be a sequence tending to infinity in the sense that for each finite subset $\Gamma_0$ of $\Gamma$, there is an $n_0 \in \sN$ such that $h(n) \not\in \Gamma_0$ if $n \ge n_0$. Clearly, if $h$ tends to infinity, then the inverse sequence $h^{-1}$ tends to infinity, too. We say that an operator $A_h \in L(l^2(\Gamma))$ is a {\em limit operator of $A \in L(l^2(\Gamma))$ defined by the sequence} $h$ if
\[
R_{h(m)}^{-1} A R_{h(m)} \to A_h \quad \mbox{and} \quad
R_{h(m)}^{-1} A^* R_{h(m)} \to A_h^*
\]
strongly as $m \to \infty$. Clearly, every operator has at most one limit operator with respect to a given sequence $h$. Note that the generating function of the shifted operator $R_r^{-1} A R_r$ is related with the generating function of $A$ by
\begin{equation} \label{equa0}
k_{R_r^{-1} A R_r} (t, \, s) = k_A (t r^{-1}, \, s r^{-1})
\end{equation}
and that the generating functions of $R_{h(m)}^{-1} A R_{h(m)}$ converge pointwise on $\Gamma \times \Gamma$ to the generating function of the limit operator $A_h$ (if the latter exists).

It is an important property of band-dominated operators that they always possess limit operators. More general, the following result can be proved by a standard Cantor diagonal argument (see \cite{RRS1,RRS2,RRSB}).
\begin{prop} \label{prop1}
Let $A$ be a band-dominated operator on $l^2(\Gamma)$. Then every sequence $h : \sN \to \Gamma$ which tends to infinity possesses a subsequence $g$ such that the limit operator $A_g$ of $A$ with respect to $g$ exists.
\end{prop}
Let $A$ be a band-dominated operator and $h : \sN \to \Gamma$ a sequence tending to infinity for which the limit operator $A_h$ of $A$ exists. Let $B$ be another band-dominated operator. By Proposition \ref{prop1} we can choose a subsequence $g$ of $h$ such that the limit operator $B_g$ exists. Then the limit operators of $A$, $A+B$ and $AB$ with respect to $g$ exist, and
\[
A_g = A_h, \qquad (A+B)_g = A_g + B_g, \qquad (AB)_g = A_g B_g.
\]
Thus, the mapping $A \mapsto A_h$ acts, at least partially, as an algebra homomorphism.

The following theorem is due to Roe \cite{Roe2}, see also \cite{RaR5}. Recall in this connection that a group $\Gamma$ is said to be {\em exact}, if its reduced translation algebra is an exact $C^*$-algebra. The latter is defined as the reduced crossed product of $l^\infty(\Gamma)$ by $\Gamma$ and coincides in our setting with the $C^*$-algebra of all band-dominated operators on $l^2(\Gamma)$. The class of exact groups is extremely rich. It includes all amenable groups (hence, all solvable groups such as the discrete Heisenberg group and the commutative groups) and all hyperbolic groups (in particular, all free groups with finitely many generators) (see \cite{Roe1}, Chapter 3).
\begin{theo}[Roe] \label{t1.1}
Let $\Gamma$ be a finitely generated discrete and exact group, and let $A$ be a band-dominated operator on $l^2(\Gamma)$. Then the operator $A$ is Fredholm on $l^2(\Gamma)$ if and only if all limit operators of $A$ are invertible and if the norms of their inverses are uniformly bounded.
\end{theo}
Note that this result holds as well if the left regular representation is replaced by the right regular one and if, thus, the operators $L_s$ and $R_t$ change their roles. In fact, in \cite{RaR5,Roe2} the results are presented in this symmetric setting. In \cite{RaR5} we showed moreover that the uniform boundedness condition in Theorem \ref{t1.1} is redundant for band operators if the group $\Gamma$ has sub-exponential growth and if not every element of $\Gamma$ is cyclic in the sense that $w^n = e$ for some positive integer $n$. For details see \cite{RaR5}. Note that the condition of sub-exponential growth is satisfied by the abelian groups $\sZ^N$, the discrete Heisenberg group and, more general, by nilpotent groups (in fact, these groups have polynomial growth), whereas the growth of the free groups $\sF_N$ is exponential.
\begin{theo} \label{t120309.1}
Let $\Gamma$ be a finitely generated discrete and exact group with sub-exponential growth which possesses at least one non-cyclic element, and let $A$ be a band operator on $l^2(\Gamma)$. Then the operator $A$ is Fredholm on $l^2(\Gamma)$ if and only if all limit operators of $A$ are invertible.
\end{theo}
\subsection{Limit operators and stability}
Let $\cY = (Y_n)$ be a sequence of finite subsets of $\Gamma$. To verify the stability of a sequence $\bA = (A_n)$ in $\cS_\cY ({\sf Sh}(\Gamma))$ via the results of the previous section, we have to choose an inflating sequence for $\cY$ and to compute the limit operators of $\op (\bA) + P_{\Gamma^\prime}$. Note that the exactness of $\Gamma$ is not relevant in this computation. Note also that large parts of this computation hold for sequences in $\cS_\cY ({\sf BDO}(\Gamma))$, too. We will consider the finite sections method for operators in $\eB \eD \eO (\Gamma)$ in detail in a forthcoming paper.

Let $\Omega$ be a finite subset of $\Gamma$ with $e \in \Omega$ which generates $\Gamma$ as a semi-group. Let $\Omega_n$ denote the set of all words with letters in $\Omega$ of length at most $n$. Thus $\Gamma = \cup_{n \ge 1} \Omega_n = \lim_{n \to 1} \Omega_n$.

By Theorem \ref{t160209.6}, the Fredholmness of the operator $\op(\bA)$ is independent of the concrete choice of the inflating sequence. For technical reasons, we choose an inflating sequence $(v_n)$ for the sequence
\[
\left( (Y_n \cup \Omega_n) (Y_n \cup \Omega_n)^{-1} (Y_n \cup \Omega_n) \right)_{n \ge 1}
\]
instead of $(Y_n)$. Since
\[
Y_n \cup \Omega_n \subset (Y_n \cup \Omega_n) (Y_n \cup \Omega_n)^{-1} \subset (Y_n \cup \Omega_n) (Y_n \cup \Omega_n)^{-1} (Y_n \cup \Omega_n),
\]
$(v_n)$ is also an inflating sequence for $(Y_n)$. Moreover, since $\lim \Omega_n = \Gamma$, one also has
\begin{equation} \label{e180209.7}
\lim \, (Y_n \cup \Omega_n) (Y_n \cup \Omega_n)^{-1} = \Gamma.
\end{equation}
Let now $\bA = (A_n) \in \cS_\cY ({\sf Sh}(\Gamma))$, set as before
\[
\op(\bA) = \sum_{n=1}^\infty R_{v_n} A_n R_{v_n}^{-1} \quad \mbox{and} \quad \Gamma^\prime = \Gamma \setminus \cup_{n=1}^\infty Y_n v_n^{-1},
\]
and let $h : \sN \to \Gamma$ be a sequence tending infinity for which the limit operator
\[
(\op(\bA) + P_{\Gamma^\prime})_h := \mbox{s-lim}_{n \to \infty} R_{h(n)}^{-1} (\op(\bA) + P_{\Gamma^\prime}) R_{h(n)}
\]
exists. Then the limit operator $(\op(\bA) + P_{\Gamma^\prime})_g$ exists for every subsequence $g$ of $h$, and it coincides with $(\op(\bA) + P_{\Gamma^\prime})_h$. So we can pass freely to subsequences of $h$. By passing to subsequences, we can restrict the computation of the limit operator to the following cases: \\[2mm]
\hspace*{5mm} {\bf Case 1:} All elements $h(n)$ belong to $\cup_{k \ge 1} \, v_k Y_k^{-1}$. \\[1mm]
\hspace*{5mm} {\bf Case 2:} No element $h(n)$ belong to $\cup_{k \ge 1} \, v_k Y_k^{-1}$. \\[2mm]
Consider {\bf Case 1}. Passing again to a subsequence of $h$ we can further suppose that each $h(n)$ belongs to one of the sets $v_k Y_k^{-1}$, say to $v_{k_n} Y_{k_n}^{-1}$, and that $v_{k_n} Y_{k_n}^{-1}$ contains no other element of the sequence $h$ besides $h(n)$. For each $n$, let $r_n$ denote the smallest non-negative integer such that $h(n) \in v_{k_n} (\partial_\Omega Y_{k_n})^{-1} \Omega_{r_n}$. Thus, $r_n$ measures the distance of $h(n)$ to the $\Omega$-boundary of $v_{k_n} Y_{k_n}^{-1}$. Finally, let $r^* := \liminf_{n \to \infty} r_n$. Again we distinguish two cases. \\[2mm]
{\bf Case 1.1: $r^*$ is finite.} Then there are infinitely many $n \in \sN$ such that $r_n = r^*$. Thus, there is a subsequence of $h$ (denoted by $h$ again) such that
\[
h(n) \in v_{k_n} Y_{k_n}^{-1} \cap v_{k_n} (\partial_\Omega Y_{k_n})^{-1} \Omega_{r^*} \quad \mbox{for all} \; n.
\]
Further, for each $n$ there is an $w_n^* \in \Omega_{r^*}$ such that $h(n) \in v_{k_n} (\partial_\Omega Y_{k_n})^{-1} w_n^*$. Since $\Omega_{r^*}$ is a finite set, one of its elements $w_n^*$ occurs for infinitely many $n$. Let $w_*$ be an element of $\Omega_{r^*}$ with this property. Consider the subsequence of $h$ which contains all elements $h(n)$ with $w_n^* = w_*$. Denoting this subsequence by $h$ again, we can hence assume that
\begin{equation} \label{e180209.1}
h(n) \in v_{k_n} Y_{k_n}^{-1} \cap v_{k_n} (\partial_\Omega Y_{k_n})^{-1} w_*
\end{equation}
for all $n$. With respect to this sequence $h$ we obtain
\begin{eqnarray} \label{e180209.3}
\lefteqn{R_{h(n)}^{-1} (\op(\bA) + P_{\Gamma^\prime}) R_{h(n)}} \nonumber \\
&& \hspace*{-3mm} = \sum_{k=1}^\infty R_{h(n)}^{-1} R_{v_k} A_k R_{v_k}^{-1} R_{h(n)} + R_{h(n)}^{-1} P_{\Gamma^\prime} R_{h(n)} \nonumber \\
&& \hspace*{-3mm} = \sum_{k \neq k_n} R_{h(n)}^{-1} R_{v_k} A_k R_{v_k}^{-1} R_{h(n)} + R_{h(n)}^{-1} P_{\Gamma^\prime} R_{h(n)} + R_{h(n)}^{-1} R_{v_{k_n}} A_{k_n} R_{v_{k_n}}^{-1} R_{h(n)}
\end{eqnarray}
with $\Gamma^\prime$ as in (\ref{e180209.2}). By (\ref{e180209.1}), $h(n) = v_{k_n} \eta_{k_n} w_*$ with $\eta_{k_n} \in (\partial_\Omega Y_{k_n})^{-1}$. Thus, the last item in (\ref{e180209.3}) becomes
\begin{equation} \label{e180209.4}
R_{w_*^{-1}} R_{\eta_{k_n}^{-1}} A_{k_n} R_{\eta_{k_n}} R_{w_*}.
\end{equation}
Set $\Pi_n := P_{(Y_{k_n} \cup \Omega_{k_n}) (Y_{k_n} \cup \Omega_{k_n})^{-1} w_*}$. By (\ref{e180209.7}), $\Pi_n \to I$ strongly.
Since $A_{k_n}$ acts on $\im P_{Y_{k_n}}$, the operator (\ref{e180209.4}) acts on $\im P_{Y_{k_n} \eta_{k_n} w_*}$. The evident inclusion
\[
Y_{k_n} \eta_{k_n} w_* \subseteq (Y_{k_n} \cup \Omega_{k_n}) (Y_{k_n} \cup \Omega_{k_n})^{-1} w_*
\]
implies that
\[
\Pi_n R_{h(n)}^{-1} R_{v_{k_n}} A_k R_{v_{k_n}}^{-1} R_{h(n)} = R_{h(n)}^{-1} R_{v_{k_n}} A_k R_{v_{k_n}}^{-1} R_{h(n)} \Pi_n =
R_{h(n)}^{-1} R_{v_{k_n}} A_k R_{v_{k_n}}^{-1} R_{h(n)}.
\]
Let now $k \neq k_n$. Then, by the inflating property,
\begin{eqnarray} \label{e180209.8}
\lefteqn{(Y_k \cup \Omega_k) (Y_k \cup \Omega_k)^{-1} (Y_k \cup \Omega_k) v_k^{-1}} \nonumber \\
&& \cap (Y_{k_n} \cup \Omega_{k_n}) (Y_{k_n} \cup \Omega_{k_n})^{-1} (Y_{k_n} \cup \Omega_{k_n}) v_{k_n}^{-1} = \emptyset.
\end{eqnarray}
Since $Y_k v_k^{-1} \subseteq (Y_k \cup \Omega_k) (Y_k \cup \Omega_k)^{-1} (Y_k \cup \Omega_k) v_k^{-1}$ and
\[
(Y_{k_n} \cup \Omega_{k_n}) (Y_{k_n} \cup \Omega_{k_n})^{-1} \eta_{k_n}^{-1} v_{k_n}^{-1} \subseteq
(Y_{k_n} \cup \Omega_{k_n}) (Y_{k_n} \cup \Omega_{k_n})^{-1} (Y_{k_n} \cup \Omega_{k_n}) v_{k_n}^{-1}
\]
we conclude from (\ref{e180209.8}) that
\[
Y_k v_k^{-1} \cap (Y_{k_n} \cup \Omega_{k_n}) (Y_{k_n} \cup \Omega_{k_n})^{-1} \eta_{k_n}^{-1} v_{k_n}^{-1} = \emptyset
\]
whence
\[
Y_k v_k^{-1} v_{k_n} \eta_{k_n} w_* \cap (Y_{k_n} \cup \Omega_{k_n}) (Y_{k_n} \cup \Omega_{k_n})^{-1} w_* = \emptyset.
\]
Since $R_{h(n)}^{-1} R_{v_k} A_k R_{v_k}^{-1} R_{h(n)}$ is an operator living on $\im P_{Y_k v_k^{-1} v_{k_n} \eta_{k_n} w_*}$, we conclude that
\[
R_{h(n)}^{-1} R_{v_k} A_k R_{v_k}^{-1} R_{h(n)} \Pi_n = \Pi_n R_{h(n)}^{-1} R_{v_k} A_k R_{v_k}^{-1} R_{h(n)} = 0
\]
for $k \neq k_n$. Hence,
\begin{eqnarray} \label{e180209.9}
\lefteqn{R_{h(n)}^{-1} (\op(\bA) + P_{\Gamma^\prime}) R_{h(n)}} \nonumber \\
&& = \sum_{k \neq k_n} R_{h(n)}^{-1} R_{v_k} A_k R_{v_k}^{-1} R_{h(n)} (I - \Pi_n) + R_{h(n)}^{-1} P_{\Gamma^\prime} R_{h(n)}
\nonumber \\
&& \qquad + \; R_{w_*}^{-1} R_{\eta_{k_n}}^{-1} A_{k_n} R_{\eta_{k_n}} R_{w_*} \Pi_n.
\end{eqnarray}
Since $\Pi_n \to I$ strongly, the first summand on the right-hand side of (\ref{e180209.9}) converges strongly (and even $^*$-strongly since $\Pi_n$ commutes with that sum) to zero. Thus,
\begin{eqnarray*}
\lefteqn{\mbox{s-lim} \, R_{h(n)}^{-1} (\op(\bA) + P_{\Gamma^\prime}) R_{h(n)}} \\
&& = \mbox{s-lim} \, R_{w_*}^{-1} R_{\eta_{k_n}}^{-1} A_{k_n} R_{\eta_{k_n}} R_{w_*} \Pi_n + \mbox{s-lim} \, R_{h(n)}^{-1} P_{\Gamma^\prime} R_{h(n)},
\end{eqnarray*}
provided that the strong limits on the right-hand side exist. The existence of the second strong limit can always be forced by passing to a suitable subsequence of $h$. Collecting these facts, we arrive at the following.
\begin{theo} \label{t180209.10}
Let $h$ be a sequence such that the limit operator $\op(\bA) + P_{\Gamma^\prime}$ exists. In Case 1.1, there is a subsequence $g$ of $h$ such that the limit operator $(P_{\Gamma^\prime})_g$ exists, and there are a monotonically increasing sequence $(k_n)$ in $\sN$, for each $n$ a vector $\eta_{k_n} \in (\partial_\Omega Y_{k_n})^{-1}$, and a $w_* \in \Gamma$ such that
\[
(\op(\bA) + P_{\Gamma^\prime})_h = \mbox{\rm s-lim} \, R_{w_*}^{-1} R_{\eta_{k_n}}^{-1} A_{k_n} R_{\eta_{k_n}} R_{w_*} + (P_{\Gamma^\prime})_g.
\]
\end{theo}
Thus, the operator $A_{k_n}$ living on $\im P_{Y_{k_n}}$ is shifted by a vector $\eta_{k_n} \in (\partial_\Omega Y_{k_n})^{-1}$ and by another vector $w_*$ independent of $n$. It is only a matter of taste to consider $A_{k_n}$ as shifted by the vector $\eta_{k_n}^{-1}$ belonging to the $\Omega$-boundary of $Y_{k_n}$. In particular, every limit operator of $\op(\bA)$ is a shift by some vector $w_*$ of a strong limit of operators
$A_{k_n}$, shifted by vectors in the boundary of $Y_{k_n}$. This is well known for the group $\sZ$ and intervals $Y_k = [-k, \, k] \cap \sZ$, and it was observed by Lindner \cite{Lin1} in case $\Gamma = \sZ^N$ and $Y_k = \Omega_k$ is a polygon with integer vertices.

Before turning to the other cases, let us specify Theorem \ref{t180209.10} to pure finite sections sequences for operators in ${\sf Sh}(\Gamma)$. The existence of the limit operator $(P_{\Gamma^\prime})_h$ is guaranteed if the strong limit
\[
\mbox{s-lim} \, R_{w_*}^{-1} R_{\eta_{k_n}}^{-1} P_{Y_{k_n}} R_{\eta_{k_n}} R_{w_*} = \mbox{s-lim} \, P_{Y_{k_n} \eta_{k_n} w_*}
\]
exists, i.e., if the set limit
\begin{equation} \label{e190209.1}
\lim Y_{k_n} \eta_{k_n} w_* =: \cY^{(h)}
\end{equation}
exists. In this case, $(P_{\Gamma^\prime})_g = I - P_{\cY^{(h)}}$.
\begin{coro}
Let $A \in {\sf Sh}(\Gamma)$, and let $h$ be a sequence such that the limit operator $\op(\bA)_h$ for the sequence $(P_{Y_n} A P_{Y_n})$ exists. In Case 1.1, there are $k_n$, $\eta_{k_n}$ and $w_*$ as in Theorem $\ref{t180209.10}$ such that the set limit
$(\ref{e190209.1})$ exists. Then
\begin{equation} \label{e190209.2}
(\op(\bA) + P_{\Gamma^\prime})_h = P_{\cY^{(h)}} A P_{\cY^{(h)}} + (I - P_{\cY^{(h)}}).
\end{equation}
Conversely, if the limit $(\ref{e190209.1})$ exists for a certain choice of $k_n$, $\eta_{k_n}$ and $w_*$ as in Theorem $\ref{t180209.10}$, then the limit operator $\op(\bA)_h$ exists for the sequence $h(n) := v_{k_n} \eta_{k_n} w_*$, and $(\ref{e190209.2})$ holds.
\end{coro}
The proof of the first assertion follows immediately from Theorem \ref{t180209.10} and from the shift invariance of the operator $A$:
\[
R_{w_*^{-1} \eta_{k_n}^{-1}} P_{Y_{k_n}} A P_{Y_{k_n}} R_{\eta_{k_n} w_*} = R_{w_*^{-1} \eta_{k_n}^{-1}} P_{Y_{k_n}} R_{\eta_{k_n} w_*} \cdot A \cdot R_{w_*^{-1} \eta_{k_n}^{-1}} P_{Y_{k_n}} R_{\eta_{k_n} w_*}.
\]
The second assertion is evident. \hfill \qed  \\[2mm]
{\bf Case 1.2: $r^*$ is infinite.} Recall that
\begin{equation} \label{e190209.3}
h(n) \in v_{k_n} Y_{k_n}^{-1} \quad \mbox{and} \quad
h(n) \not\in v_{k_n} (\partial_\Omega Y_{k_n})^{-1} \Omega_{r_n -1}
\end{equation}
for all $n \in \sN$. The second assertion in (\ref{e190209.3}) implies that
\[
h(n) \Omega_{r_n -1}^{-1} \cap v_{k_n} (\partial_\Omega Y_{k_n})^{-1} = \emptyset.
\]
Hence, we can rewrite (\ref{e190209.3}) as
\begin{equation} \label{e190209.4}
e \in Y_{k_n} v_{k_n}^{-1} h(n) \quad \mbox{and} \quad
\Omega_{r_n -1} \cap (\partial_\Omega Y_{k_n}) v_{k_n}^{-1} h(n) = \emptyset.
\end{equation}
We claim that this implies that
\begin{equation} \label{e190209.5}
\Omega_{r_n -1} \subseteq Y_{k_n} v_{k_n}^{-1} h(n).
\end{equation}
Suppose (\ref{e190209.5}) is wrong. Then $\Omega_{r_n -1}$ has at least one point outside $Y_{k_n} v_{k_n}^{-1} h(n)$, say $a$, but it also has points inside this set, for example the point $e$ due to the first assumption of (\ref{e190209.4}). Write $a$ as a product $a = w_{r_n-1} \ldots w_1 w_0$ of elements $w_i \in \Omega$ with $w_0 := e$, and let $0 \le j < r_n-1$ be the smallest integer such that
\[
w_j \ldots w_1 w_0 \in Y_{k_n} v_{k_n}^{-1} h(n), \quad \mbox{but} \quad  w_{j+1} w_j \ldots w_1 w_0 \not\in Y_{k_n} v_{k_n}^{-1} h(n).
\]
Then $\Omega w_j \ldots w_1 w_0 \not\subseteq Y_{k_n} v_{k_n}^{-1} h(n)$, hence
\[
w_j \ldots w_1 w_0 \in \partial_\Omega (Y_{k_n} v_{k_n}^{-1} h(n)).
\]
Since $w_j \ldots w_1 w_0 \in \Omega_{r_n -1}$, this contradicts
the second assertion of (\ref{e190209.4}), and the claim
(\ref{e190209.5}) follows. Roughly speaking, we used the fact that $\Omega$-boundaries do not have gaps. Since $P_{\Omega_n} \to I$ strongly, we conclude from (\ref{e190209.5}) that
\begin{equation} \label{e190209.8}
P_{Y_{k_n} v_{k_n}^{-1} h(n)} \to I \quad \mbox{strongly}.
\end{equation}
\begin{theo} \label{t190209.6}
Let $\bA \in \cS_\cY ({\sf Sh}(\Gamma))$, and let $h$ be a sequence such that the limit operator $\op(\bA)_h$ exists. Then in Case 1.2,
\begin{equation} \label{e190209.7}
\op(\bA)_h = A \quad \mbox{with} \quad A := \mbox{\rm s-lim} A_n P_{Y_n}.
\end{equation}
\end{theo}
{\bf Proof.} It is sufficient to prove (\ref{e190209.7}) for pure finite sections sequences $\bA = (P_{Y_n} A P_{Y_n})$ with $A \in {\sf Sh}(\Gamma)$. For these sequences, one has
\begin{eqnarray*}
R_{h(n)}^{-1} (\op(\bA) + P_{\Gamma^\prime}) R_{h(n)}
& = & \sum_{k \neq k_n} R_{h(n)}^{-1} R_{v_k} P_{Y_k} A P_{Y_k} R_{v_k}^{-1} R_{h(n)} (I - P_{Y_{k_n} v_{k_n}^{-1} h(n)}) \\
&& \quad + \; R_{h(n)}^{-1} P_{\Gamma^\prime} R_{h(n)} (I - P_{Y_{k_n} v_{k_n}^{-1} h(n)}) \\
&& \quad + \; P_{Y_{k_n} v_{k_n}^{-1} h(n)} A P_{Y_{k_n} v_{k_n}^{-1} h(n)}.
\end{eqnarray*}
Letting $n$ go to infinity the assertion follows due to (\ref{e190209.8}). \hfill \qed \\[3mm]
Thus, in Case 1.2, the invertibility of the limit operators of $\op(\bA) + P_{\Gamma^\prime}$

 follows already from the invertibility of $A$. \\[2mm]
Now consider {\bf Case 2}, i.e., suppose that none of the $h(n)$ belongs to $\cup v_k Y_k^{-1}$. For $n \in \sN$, let $r_n$ stand for the smallest non-negative integer such that there is a $k_n \in \sN$ with $h(n) \in v_{k_n} (\partial_\Omega Y_{k_n})^{-1} \Omega_{r_n}$. Consequently,
\[
h(n) \not\in v_{k_n} (\partial_\Omega Y_{k_n})^{-1} \Omega_{r_n - 1} \quad \mbox{for all} \; n.
\]
Again we set $r^* := \liminf r_n$ and distinguish two cases. \\[2mm]
{\bf Case 2.1: $r^*$ is finite.} We proceed as in Case 1.1 and find a subsequence of $h$ (denoted by $h$ again) and an element $w_* \in \Gamma$ such that
$h(n) \in v_{k_n} (\partial_\Omega Y_{k_n})^{-1} w_*$. Since the inclusion $h(n) \in v_{k_n} Y_{k_n}^{-1}$ in (\ref{e180209.1}) had not been used in Case 1.1 we can continue exactly as in that case to obtain that Theorem \ref{t180209.10} and its corollary hold verbatim in the case at hand, too. \\[2mm]
{\bf Case 2.2: $r^*$ is infinite.} As in Case 1.2, we choose the sequence $(r_n)$ as strongly monotonically increasing. Then we have
\begin{equation} \label{e190209.9}
h(n) \not\in v_k Y_k^{-1} \quad \mbox{for all} \; k, \, n,
\end{equation}
\begin{equation} \label{e190209.10}
h(n) \not\in v_k (\partial_\Omega Y_k)^{-1} \Omega_{r_n -1} \quad \mbox{for all} \; k, \, n.
\end{equation}
We claim that these two facts imply that
\begin{equation} \label{e190209.12}
\Omega_{r_n -1} \cap Y_k v_k^{-1} h(n) = \emptyset \quad \mbox{for all} \; k, \, n.
\end{equation}
Indeed, from (\ref{e190209.9}) we conclude that $e \not\in Y_k v_k^{-1} h(n)$. Thus, for each $k$ and $n$, $\Omega_{r_n -1}$ contains points from the complement of $Y_k v_k^{-1} h(n)$, for instance the point $e$. Suppose that $\Omega_{r_n -1}$ also contains points in $Y_k v_k^{-1} h(n)$. Then the arguments from Case 1.2 imply that $\Omega_{r_n -1}$ contains points in the $\Omega$-boundary of $Y_k v_k^{-1} h(n)$. But (\ref{e190209.10}) implies that
$\Omega_{r_n -1} \cap (\partial_\Omega Y_k) v_k^{-1} h(n) = \emptyset$. Thus, $\Omega_{r_n -1}$ is completely located in the  complement of $Y_k v_k^{-1} h(n)$, whence (\ref{e190209.12}).

Since the operator $R_{h(n)}^{-1} R_{v_k} A_k R_{v_k}^{-1} R_{h(n)}$ lives on $\im P_{Y_k v_k^{-1} h(n)}$, we obtain from (\ref{e190209.12})
\begin{eqnarray*}
R_{h(n)}^{-1} (\op(\bA) + P_{\Gamma^\prime}) R_{h(n)}
& = & \sum_{k \ge 1} R_{h(n)}^{-1} R_{v_k} A_k R_{v_k}^{-1} R_{h(n)} (I - P_{\Omega_{r_n - 1}}) \\
&& \quad + \; R_{h(n)}^{-1} P_{\Gamma^\prime} R_{h(n)} (I - P_{\Omega_{r_n - 1}}) + P_{\Omega_{r_n - 1}}.
\end{eqnarray*}
The first two summands on the right-hand side of this equality tend strongly to zero as $n \to \infty$, whereas the third one tends strongly to the identity. Thus, the identity operator is the only limit operator of $\op(\bA) + P_{\Gamma^\prime}$ in Case 2.2. The following theorem summarizes the results from Cases 1.1 - 2.2. %
\begin{theo} \label{t190209.13}
Let $\bA \in \cS_\cY ({\sf Sh}(\Gamma))$. Then the limit operators of $\op(\bA) + P_{\Gamma^\prime}$ are the identity operator $I$,
the operator $A := \mbox{\rm s-lim} \, A_n P_{Y_n}$, and all operators of the form
\[
\mbox{\rm s-lim} \, R_{w_*}^{-1} R_{\eta_{k_n}}^{-1} A_{k_n}
R_{\eta_{k_n}} R_{w_*} + (P_{\Gamma^\prime})_g
\]
with a suitable subsequence $g$ of $h$ and with elements $\eta_{k_n} \in (\partial_\Omega Y_{k_n})^{-1}$ and $w_* \in \Gamma$.
\end{theo}
Combining this theorem with Theorems \ref{t160209.6} $(b)$ and \ref{t1.1} we arrive at the following stability results.
\begin{theo} \label{t190209.14}
Let $\Gamma$ be an exact discrete group, and let $\bA \in \cS_\cY ({\sf Sh}(\Gamma))$. The sequence $\bA$ is stable if and only if the operator $A := \mbox{\rm s-lim} \, A_n P_{Y_n}$ and all operators of the form
\[
\mbox{\rm s-lim} \,  R_{\eta_{k_n}}^{-1} A_{k_n} R_{\eta_{k_n}} + R_{w_*} (P_{\Gamma^\prime})_g R_{w_*}^{-1}
\]
with a suitable subsequence $g$ of $h$ and with elements $\eta_{k_n} \in (\partial_\Omega Y_{k_n})^{-1}$ and $w_* \in \Gamma$ are invertible and if the norms of their inverses are uniformly bounded.
\end{theo}
\begin{coro} \label{c190209.15}
Let $\Gamma$ be an exact discrete group, and let $A \in {\sf Sh}(\Gamma)$. The sequence $\bA = (P_{Y_n} A P_{Y_n})$ is stable if and only if the operator $A$ and all operators
\[
P_{\cY^{(h)}} A P_{\cY^{(h)}} : \im P_{\cY^{(h)}} \to \im P_{\cY^{(h)}}
\]
where
\begin{equation} \label{e190209.16}
\cY^{(h)} := \lim Y_{k_n} \eta_{k_n}
\end{equation}
with certain elements $\eta_{k_n} \in (\partial_\Omega Y_{k_n})^{-1}$ are invertible and if the norms of their inverses are uniformly bounded.
\end{coro}
Theorem \ref{t120309.1} allows us to remove the uniform boundedness condition in the previous corollary.
\begin{coro} \label{c120210.1}
Let $\Gamma$ be a finitely generated discrete and exact group with sub-exponential growth which possesses at least one non-cyclic element, and let $A \in {\sf Sh}(\Gamma)$ be a band operator. Then the sequence $\bA = (P_{Y_n} A P_{Y_n})$ is stable if and only if the operators mentioned in the previous corollary are invertible.
\end{coro}
\subsection{Geodesic paths}
Now we turn to special sequences $\cY = (Y_n)$ and $\eta : \sN \to \Gamma$ for which the existence of the set limit (\ref{e190209.16}) can be guaranteed. Let again $\Omega_n$ refer to the set of all products of at most $n$ elements of $\Omega$ and set $\Omega_0 := \{e\}$. A sequence $(\nu_n)$ in $\Gamma$ is called a {\em geodesic path} (with respect to $\Omega$) if there is a sequence $(w_n)$ in $\Omega \setminus \{e\}$ such that $\nu_n = w_1 w_2 \ldots w_n$ and $\nu_n \in \Omega_n \setminus \Omega_{n-1}$ for each $n \ge 1$. Note that this condition implies that each $\nu_n$ is in the {\em right $\Omega$-boundary of} $\Omega_n$, which is the set of all $w \in \Omega_n$ for which $w \Omega$ is not a subset of $\Omega_n$.

We will see now that the $\lim \Omega_n \eta_n$ exists if $\eta$ is an inverse geodesic path, i.e., if $\eta_n = \nu_n^{-1}$ for a geodesic path $\nu$.
\begin{lemma} \label{l190209.17}
Let $(w_n)_{n \ge 1}$ be a sequence in $\Omega$ and set $\eta_n := w_n^{-1} w_{n-1}^{-1} \ldots w_1^{-1}$ for $n \ge 1$. Then the set limit $\lim \Omega_n \eta_n$ exists, and
\begin{equation} \label{e190209.18}
\lim \Omega_n \eta_n = \cup_{n \ge 1} \Omega_n \eta_n.
\end{equation}
\end{lemma}
{\bf Proof.} For $n \ge 1$, one has $\Omega_n \eta_n = \Omega_n w_{n+1} w_{n+1}^{-1} w_n^{-1} \ldots w_1^{-1} \subseteq \Omega_{n+1} \eta_{n+1}$.
These inclusions imply the existence of the set limit and the equality (\ref{e190209.18}). \hfill \qed \\[3mm]
The natural question arises whether every sequence $\eta : \sN \to \Gamma$ for which the set limit (\ref{e190209.16}) exists has a subsequence which is a subsequence of an inverse geodesic path. If the answer is affirmative, then it would prove sufficient to consider strong limits with respect to inverse geodesic paths in Theorem \ref{t190209.14} and its corollary. We are going to answer this question for two special families of groups.
\subsection{Commutative groups}
Let $\Gamma$ be a commutative group which is generated, as a semi-group, by the finite set $\Omega$ with $e \in \Omega$.
Define $\Omega_n$ as in the previous section.
\begin{prop} \label{p190209.19}
Let $\Gamma$ be commutative, and let $\mu = (\mu_n)_{n \in \sN}$ be a sequence in $\Gamma$ which has a subsequence $(\mu_n)_{n \in \sN_0}$ with $\mu_n \in \Omega_n \setminus \Omega_{n-1}$ for each $n \in \sN_0$. Then $(\mu_n)_{n \in \sN_0}$ has a subsequence which is a subsequence of a geodesic path.
\end{prop}
{\bf Proof.} Let $\Omega = \{ e, \, \omega_1, \, \ldots, \, \omega_k \}$. Each $\mu_n$ can be written as $\omega_1^{e_{1n}} \omega_2^{e_{2n}} \ldots \omega_k^{e_{kn}}$ where $e_{1n} + e_{2n} + \ldots + e_{kn} = n$ for $n \in \sN_0$. (We do not claim that this representation of $\mu_n$ is unique.) Consider the sequence $(e_{1n})_{n \in \sN_0}$. This sequence has a constant subsequence or a strongly monotonically increasing subsequence. Let $(e_{1n})_{n \in \sN_1}$ with an infinite subset $\sN_1$ of $\sN_0$ be a subsequence of $(e_{1n})_{n \in \sN_0}$ which owns one of these properties. Then consider $(e_{2n})_{n \in \sN_1}$ and choose a subsequence $(e_{2n})_{n \in \sN_2}$ which is constant or strongly monotonically increasing. We proceed in this way. After $k$ steps we arrive at a subsequence $(\mu_n)_{n \in \sN_k}$ of $(\mu_n)_{n \in \sN_0}$ with $\mu_n = \omega_1^{f_{1n}} \omega_2^{f_{2n}} \ldots \omega_k^{f_{kn}}$ and $f_{1n} + f_{2n} + \ldots + f_{kn} = n$ for $n \in \sN_k$ and where each of the sequences $(f_{in})$ is either constant or strongly monotonically increasing.

For $n \in \sN_k$ let $\nu_n := \mu_n$, and set $\nu_0 := e$. Let $(k_n)$ be the enumeration of the elements of $\sN_k$ in increasing order, and set $k_0 := 0$. In order to define $\nu_n$ for $k_r < n < k_{r+1}$ we proceed as follows. Let $i_1$ be the smallest positive integer such that $f_{i_1 k_r} < f_{i_1  k_{r+1}}$. For $l = 1, \, \ldots, \, f_{i_1 k_{r+1}} - f_{i_1 k_r}$, set
\[
\nu_{k_r + l} := \omega_1^{f_{1 k_r}} \ldots \, \omega_{i_1-1}^{f_{i_1-1, k_r}} \omega_{i_1}^{f_{i_1, k_r + l}} \omega_{i_1+1}^{f_{i_1+1, k_r}} \ldots \, \omega_k^{f_{k k_r}}.
\]
Now we are looking for the next subscript, say $i_2$, for which
the exponents at $\omega_{i_2}$ of $\nu_{k_n}$ and $\nu_{k_{n+1}}$ are different and proceed in the same way. After a finite number of steps, we arrive at a sequence $\nu = (\nu_n)_{n \in \sN}$
with $\nu_n \in \Omega_n$ for each $n \in \sN$.

It remains to show that the sequence $\nu$ is a geodesic path, i.e. that $\nu_n \in \Omega_n \setminus \Omega_{n-1}$ for each $n$. Suppose that $\nu_k \not\in \Omega_k \setminus \Omega_{k-1}$ for some $k \ge 2$. Then $\nu_k$ is a product of $l < k$ elements from $\Omega \setminus \{e\}$. Choose $n$ such that $k_n > k$ and let $a \in \Omega_{k_n - k}$ such that $\mu_{k_n} = a \nu_k$. Then
\[
\mu_{k_n} \in \Omega_{k_n - k} \Omega_l = \Omega_{k_n - k + l}
\quad \mbox{with} \quad k_n - k + l < k_n,
\]
a contradiction to the hypothesis that $\mu_{k_n} \in \Omega_{k_n} \setminus \Omega_{k_n -1}$ for each $n$. \hfill \qed \\[3mm]
Since commutative groups are exact, one has the following consequences.
\begin{coro} \label{c220209.1}
Let $\Gamma$ be a commutative discrete group, and let $\Omega$ be a finite subset of $\Gamma$ which generates $\Gamma$ as a semi-group. Set $Y_n := \Omega_n$, and let $\bA \in \cS_\cY ({\sf Sh}(\Gamma))$. The sequence $\bA$ is stable if and only if the operator $A := \mbox{\rm s-lim} \, A_n P_{\Omega_n}$ and, for each inverse geodesic path $\eta$, the operator
\[
\mbox{\rm s-lim} \,  R_{\eta_n}^{-1} A_n R_{\eta_n} : \im P_{\cup \Omega_n \eta_n} \to \im P_{\cup \Omega_n \eta_n}
\]
are invertible and if the norms of the inverses of these operators are uniformly bounded.
\end{coro}
\begin{coro} \label{c220209.2}
Let $\Gamma$ and $\Omega$ be as in Corollary $\ref{c220209.1}$, and let $A \in {\sf Sh}(\Gamma)$. The sequence $\bA = (P_{\Omega_n} A P_{\Omega_n})$ is stable if and only if the operator $A$ and, for each inverse geodesic path $\eta$, the operator
\[
P_{\cup \Omega_n \eta_n} A P_{\cup \Omega_n \eta_n} : \im P_{\cup \Omega_n \eta_n} \to \im P_{\cup \Omega_n \eta_n}
\]
are invertible and if the norms of their inverses are uniformly bounded.
\end{coro}
In many cases, there will be only finitely many different set limits $\lim \Omega_n \eta_n$; then the uniform boundedness condition in the previous corollaries is redundant. The same happens if $A$ is a band operator by Theorem \ref{t120309.1}.

The perhaps most important consequence of Corollary $\ref{c220209.1}$ is that the finite sections method for operators in ${\sf Sh}(\Gamma)$ is fractal. More general, one has the following.
\begin{coro} \label{c220209.3}
Let $\Gamma$, $\Omega$ and $\cY$ be as in Corollary $\ref{c220209.1}$. Then the algebra $\cS_\cY ({\sf Sh}(\Gamma))$ is fractal.
\end{coro}
Roughly saying, an algebra of matrix sequences is fractal if each sequence in the algebra can be reconstructed from each of its (infinite) subsequences modulo a sequence tending to zero in the norm. For an exact definition and some properties of sequences in fractal algebras, see \cite{HRS2,Roc0}. The proof of Corollary \ref{c220209.3} follows immediately from Corollary \ref{c220209.1}. See Theorem 1.69 in \cite{HRS2} and its corollary for the argument.
\subsection{The free non-commutative group $\sF_N$}
Proposition \ref{p190209.19} does certainly not hold for all discrete groups. For example, let $\Gamma = \sF_2$ with generators $u$ and $v$, set $\Omega := \{e, \, u^{\pm 1}, \, v^{\pm 1} \}$, and let $\Omega_n$ stand for the set of all products of at most $n$ elements of $\Omega$. Consider $\eta_n := vu^{n-1}$. It is easy to see (indeed, drawing pictures will help a lot in what follows) that the set limit $\lim \Omega_n \eta_n$ exists, but the sequence $\eta$ has no subsequence which is a subsequence of an inverse geodesic path. On the other hand, a simple calculation gives
\[
\lim \Omega_n \eta_n = \lim \Omega_{n-1} u^{n-1};
\]
thus, the set limit $\lim \Omega_n \eta_n$ coincides with another set limit which {\em is} taken with respect to an inverse geodesic path. We will see now that this observation is archetypal for the free non-commutative groups $\sF_N$.

Still for a moment, let $\Gamma$ be a general discrete group with a finite set $\Omega$ of generators. Let $(\eta_{k_n})$ be a sequence with
$\eta_{k_n} \in (\Omega_{k_n} \setminus \Omega_{k_n - 1})^{-1}$ for each $n$. Write $\eta_{k_n}^{-1}$ as
\begin{equation} \label{e020309.1}
\eta_{k_n}^{-1} = \omega_1^{(n)} \omega_2^{(n)} \ldots \, \omega_{k_n}^{(n)} \quad \mbox{with} \quad \omega_i^{(n)} \in \Omega \setminus \{e\}
\end{equation}
for each $i = 1, \, \ldots, \, k_n$. Again, we do not claim that this representation is unique. Since $\Omega$ is finite, there is an $\tilde{\omega}_1 \in \Omega$ such that $\omega_1^{(n)} = \tilde{\omega}_1$ for infinitely many $n \in \sN$, say for all $n \in \sN_1$. By the same argument, there is an $\tilde{\omega}_2 \in \Omega$ such that $\omega_2^{(n)} = \tilde{\omega}_2$ for infinitely many $n \in \sN_1$, say for all $n \in \sN_2$. We proceed in that way to obtain a sequence $(\tilde{\omega}_n)_{n \in \sN}$ in $\Omega \setminus \{e\}$ having the property that, for each $r \in \sN$, there are infinitely many elements $\eta_{k_n}$ with
\begin{equation} \label{e020309.2}
\eta_{k_n}^{-1} = \tilde{\omega}_1 \tilde{\omega}_2 \ldots \, \tilde{\omega}_r \, \omega_{r+1}^{(n)} \ldots \, \omega_{k_n}^{(n)}.
\end{equation}
For $r \ge 1$, set
\begin{equation} \label{e020309.3}
\tilde{\eta}_r := (\tilde{\omega}_1 \tilde{\omega}_2 \ldots \, \tilde{\omega}_r)^{-1}.
\end{equation}
By Lemma \ref{l190209.17}, the set limit $\lim \Omega_r \tilde{\eta}_r$ exists.
\begin{lemma} \label{l020309.4}
Let $\eta_{k_n}$ and $\tilde{\eta}_r$ be as in $(\ref{e020309.1})$ and $(\ref{e020309.3})$, respectively. Then
\begin{equation} \label{e020309.5}
\lim_{r \to \infty} \Omega_r \tilde{\eta}_r \subseteq \limsup_{n \to \infty} \Omega_{k_n} \eta_{k_n}.
\end{equation}
\end{lemma}
{\bf Proof.} Let $x \in \Omega_r \tilde{\eta}_r$ for some $r$, and let $\eta_{k_n}^{-1}$ be as in (\ref{e020309.2}). Then
\[
x \in \Omega_r \tilde{\eta}_r =  \Omega_r \omega_{r+1}^{(n)} \ldots \, \omega_{k_n}^{(n)} \, (\omega_{k_n}^{(n)})^{-1} \ldots \, (\omega_{r+1}^{(n)})^{-1} \tilde{\omega}_r^{-1} \ldots \, \tilde{\omega}_1^{-1} \subseteq \Omega_{k_n} \eta_{k_n}.
\]
Since there are infinitely many elements as in (\ref{e020309.2}), this inclusion implies that $x \in \limsup \Omega_{k_n} \eta_{k_n}$, whence
$\cup_{r \ge 1} \Omega_r \tilde{\eta}_r \subseteq \limsup \Omega_{k_n} \eta_{k_n}$. This is the assertion. \hfill \qed \\[3mm]
It is one consequence of the lemma that the set limits $\lim \Omega_{k_n} \eta_{k_n}$ cannot be too small. In particular, they contain a shifted copy of $\Omega_r$ for each $r$ and are, thus, growing sets in the sense of Shteinberg (see \cite{Ste1} and Definition 2.4.8 in \cite{RRSB}).

In general, one cannot expect that equality holds in (\ref{e020309.5}). For example, let $\Gamma$ be the (additively written) group $\sZ^2$ with $\Omega = \{(0, \, 0), \, (\pm 1, \, 0), \, (0, \, \pm 1)\}$ and consider the sequence $\eta_{2n} = (-n, \, -n)$. If we write $- \eta_{2n}$ as
\[
- \eta_{2n} = (1, \, 0) + \ldots + (1, \, 0) + (0, \, 1)  + \ldots + (0, \, 1)
\]
with each summand occurring $n$ times, then the above construction yields $\tilde{\eta}_r := (-r, \, 0)$. In this setting, both set limits $\lim \Omega_{2n} \eta_{2n}$ and $\lim \Omega_r \tilde{\eta}_r$ exist, but they do not coincide (the first one is the intersection of $\sZ^2$ with a half plane, the second one with a quadrant).

It turns out that, in case of the free non-commutative groups $\sF_N$, equality holds in (\ref{e020309.5}).
\begin{theo} \label{t020309.6}
For $N > 1$, let $\sF_N$ be the free group generated by its elements $\omega_1, \, \ldots, \, \omega_N$, set $\Omega := \{e, \, \omega_1^{\pm 1}, \, \ldots, \, \omega_N^{\pm 1}\}$, and let $\Omega_n$ be the set of all products of elements of $\Omega$ of length at most $n$. Further let $(\eta_{k_n})$ be a sequence with
\[
\eta_{k_n} \in (\Omega_{k_n} \setminus \Omega_{k_n - 1})^{-1}
\]
which we write as in $(\ref{e020309.2})$ and let $(\tilde{\eta}_r)$ be the associated sequence as in $(\ref{e020309.3})$. Then
\begin{equation} \label{e020309.7}
\liminf_{n \to \infty} \Omega_{k_n} \eta_{k_n} \subseteq \lim_{r \to \infty} \Omega_r \tilde{\eta}_r.
\end{equation}
In particular, if the set limit $\lim_{n \to \infty} \Omega_{k_n} \eta_{k_n}$ exists, then
\begin{equation} \label{e020309.8}
\lim_{n \to \infty} \Omega_{k_n} \eta_{k_n} = \lim_{r \to \infty} \Omega_r \tilde{\eta}_r.
\end{equation}
\end{theo}
{\bf Proof.} Let $x \in \liminf \Omega_{k_n} \eta_{k_n}$. Then there is an $n_0$ such that $x \in \Omega_{k_n} \eta_{k_n}$ for $n \ge n_0$. Thus, for each $n \ge n_0$,
\[
x \in \Omega_{k_n} (\omega_{k_n}^{(n)})^{-1} \ldots \, (\omega_{n+1}^{(n)})^{-1} \, \tilde{\omega}_n^{-1} \ldots \,
\tilde{\omega}_1^{-1}.
\]
Choose elements $\nu_i^{(n)}$ in $\Omega$ such that
\begin{equation} \label{e020309.9}
x = \underbrace{\nu_1^{(n)} \ldots \, \nu_{k_n}^{(n)}}_{(**)} \; \underbrace{ (\omega_{k_n}^{(n)})^{-1} \ldots \, (\omega_{n+1}^{(n)})^{-1} \, \tilde{\omega}_n^{-1} \ldots \,
\tilde{\omega}_1^{-1}}_{(*)}.
\end{equation}
The assumption $\eta_{k_n} \in (\Omega_{k_n} \setminus \Omega_{k_n - 1})^{-1}$ guarantees that there is no cancelation possible inside part $(*)$ of the representation (\ref{e020309.9}) but, of course, there might be cancelation inside part $(**)$ as well as between the most right of the $\nu$ and the most left of the $\omega^{-1}$.

For each $n \ge n_0$, we cancel the representation (\ref{e020309.9}) of $x$ as far as possible. Suppose that, after complete cancelation, at least one factor $(\omega_k^{(n)})^{-1}$ remains in each representation. Then, for each $n \ge 1$, we can represent $x$ as a word without further cancelation which starts
from the right-hand side with $\ldots \tilde{\omega}_n^{-1} \ldots \, \tilde{\omega}_1^{-1}$ and, hence, has length at most $n$. This is impossible since each $x \in \sF_N$ can be uniquely represented as a reduced word of finite length. This contradiction shows that there is at least one $n \ge n_0$ such that all factors $(\omega_k^{(n)})^{-1}$ in the representation (\ref{e020309.9}) can be canceled. Thus, $x \in \Omega_k \, \tilde{\omega}_k^{-1} \ldots \, \tilde{\omega}_1^{-1} = \Omega_k \, \tilde{\eta}_k$ for some $k \ge n_0$. Since the set sequence $(\Omega_k \, \tilde{\eta}_k)$ is monotonically increasing, this implies
\[
x \in \cup_{k \ge 1} \, \Omega_k \, \tilde{\eta}_k = \lim \Omega_k \, \tilde{\eta}_k
\]
whence the first assertion. Combining this result with Lemma \ref{l020309.4}, the second assertion follows. \hfill \qed \\[3mm]
Thus, each set limit $\lim \Omega_{k_n} \eta_{k_n}$ can be obtained as a set limit along an inverse geodesic path. Since free groups are exact, this leads to the same consequences as for commutative groups.
\begin{coro} \label{c020309.10}
Let $\Gamma = \sF_N$ and $\Omega$ and $\Omega_n$ as in Theorem $\ref{t020309.6}$. Set $Y_n := \Omega_n$, and let $(A_n) \in \cS_\cY ({\sf Sh}(\sF_N))$. The sequence $(A_n)$ is stable if and only if the operator $A := \mbox{\rm s-lim} \, A_n P_{\Omega_n}$ and, for each inverse geodesic path $\eta$, the
operator
\[
\mbox{\rm s-lim} \,  R_{\eta_n}^{-1} A_n R_{\eta_n} : \im P_{\cup \Omega_n \eta_n} \to \im P_{\cup \Omega_n \eta_n}
\]
are invertible and if the norms of the inverses of these operators are uniformly bounded.
\end{coro}
\begin{coro} \label{c020309.11}
Let $\Omega$ be as in Corollary $\ref{c020309.10}$ and let $A \in {\sf Sh}(\sF_N)$. The sequence $\bA = (P_{\Omega_n} A P_{\Omega_n})$ is stable if and only if the operators $A$ and, for each inverse geodesic path $\eta$,
\[
P_{\cup \Omega_n \eta_n} A P_{\cup \Omega_n \eta_n} : \im P_{\cup \Omega_n \eta_n} \to \im P_{\cup \Omega_n \eta_n}
\]
are invertible and if the norms of their inverses are uniformly bounded.
\end{coro}
\begin{coro} \label{c020309.13}
Let $\Omega, \, \cY$ be as in Corollary $\ref{c020309.10}$. Then the algebra $\cS_\cY ({\sf Sh}(\sF_N))$ is fractal.
\end{coro}
{\small Author's address: \\[3mm]
Steffen Roch, Technische Universit\"at Darmstadt, Fachbereich
Mathematik, Schlossgartenstrasse 7, 64289 Darmstadt,
Germany. \\
E-mail: roch@mathematik.tu-darmstadt.de}
\end{document}